\documentclass[letterpaper,11pt,onecolumn]{article}
\usepackage[margin=0.8in]{geometry}

\usepackage[T1]{fontenc}
\usepackage[utf8]{inputenc}
\usepackage[english]{babel}
\usepackage{amsmath}
\usepackage{amsbsy}
\usepackage{amsfonts}
\usepackage{amssymb}
\usepackage{array}
\usepackage{bm}
\usepackage{booktabs}
\usepackage{csquotes}
\usepackage{enumerate}
\usepackage{float}
\usepackage{graphicx}
\usepackage{placeins}
\usepackage{mathrsfs}
\usepackage{pgf}
\usepackage{pgfplots}
\usepackage{color}
\usepackage{tikz}
\usepackage[font={footnotesize}]{caption}
\usepackage{subcaption}

\usepackage[pdftex]{hyperref}
\hypersetup{colorlinks=false}
\usepackage{cleveref}
\crefname{figure}{fig}{figures}

\def\be{\begin{equation}}
\def\ee{\end{equation}}
\def\ba{\begin{array}}
\def\ea{\end{array}}
\def\bea{\begin{eqnarray}}
\def\eea{\end{eqnarray}}
\def\beas{\begin{eqnarray*}}
\def\eeas{\end{eqnarray*}}

\newcommand{\half}{\frac{1}{2}}

\DeclareUnicodeCharacter{FFFD}{\textquestiondown}

\newcommand{\AB}{\bm{A}}

\newcommand{\ptl}{{\partial}}

\newcommand{\ds}{\displaystyle}

\catcode `@=11

\@definecounter{Remark}

\newbox \itemlist@label
\newdimen \itemlist@labelpad

\def\itemlist@makelabel#1{%
\setbox\itemlist@label =\hbox{#1}%
\ifdim \wd\itemlist@label >\labelwidth
\itemlist@labelpad=\textwidth
\advance\itemlist@labelpad by -\rightmargin
\advance\itemlist@labelpad by -\@totalleftmargin
\advance\itemlist@labelpad by \labelwidth
\hbox to \itemlist@labelpad {#1\hfill}%
\else #1\hfill
\fi
}

\let\div\relax
\DeclareMathOperator{\div}{div}

\makeindex

\numberwithin{equation}{section}

\graphicspath{{figures/}}

\usepackage[backend=bibtex,
	giveninits=true,
	maxbibnames=4,
	date=year,
	doi=false,
	url=false,
	isbn=false,
	eprint=false]{biblatex}

\title{A Space-Time Discontinuous Petrov-Galerkin Finite Element Formulation for a Modified Schr\"{o}dinger Equation for Laser Pulse Propagation in Waveguides}
\date{\vspace{-5ex}}
\begin{document}
\maketitle
\begin{center}
{{\bf  
{
A. Chakraborty$^a$ \footnote{Corresponding author: ankit.chakraborty@austin.utexas.edu, ORCID:0  },
J.  Mu\~{n}oz-Matute$^{a,b}$,\\[8pt]
L. Demkowicz$^a$ 
and J. Grosek$^c$
\\[12pt]}}}
{
$^a$Oden Institute, The University of Texas at Austin, USA,
$^b$BCAM,
$^c$Air Force Research Laboratory
}
\end{center}

\begin{abstract}
In this article, we propose a modified nonlinear Schr\"{o}dinger equation for modeling pulse propagation in optical waveguides. The proposed model bifurcates into a system of elliptic and hyperbolic equations depending on waveguide parameters. The proposed model leads to a stable first-order system of equations, distinguishing itself from the canonical nonlinear Schr\"{o}dinger equation. We have employed the space-time discontinuous Petrov-Galerkin finite element method to discretize the first-order system of equations. We present a stability analysis for both the elliptic and hyperbolic systems of equations and demonstrate the stability of the proposed model through several numerical examples on space-time meshes.
\end{abstract}



\subsection*{Acknowledgments}  A. Chakraborty and L. Demkowicz were supported with AFOSR grant FA9550-19-1-0237, J. Mu\~{n}oz-Matute has received funding from the European Union's Horizon 2020 research and innovation programme under the Marie Sklodowska-Curie individual fellowship grant agreement No. 101017984 (GEODPG). 


\section{Introduction}
Pulse propagation in optical waveguides is a fundamental topic in photonics, with significant implications for telecommunications \cite{photonics_1,photonics_2,photonics_3}, medical imaging \cite{photonics_4,photonics_5,photonics_2}, and various other applications \cite{photonics_7,photonics_8}. The nonlinear Schr\"{o}dinger (NLS) equation has long been a cornerstone in modeling the dynamics of pulse propagation \cite{NLS_1,NLS_2,NLS_3,NLS_4} in optical waveguides. It captures features of pulse evolution, including dispersion and nonlinear effects, making it indispensable for both theoretical and numerical investigations. The NLS equation is particularly well-suited for describing the propagation of optical pulses in fibers, where the balance between dispersive spreading and nonlinear self-focusing leads to the formation of stable structures known as solitons. These solitons maintain their shape over long distances, which is crucial for applications in long-haul optical communication systems \cite{NLS_LH_1,NLS_LH_2,NLS_4}. \\
In terms of computational methods employed to solve the NLS equation, traditional approaches such as the split-step Fourier method \cite{SS_2, SS_3, SS_4, SS_1} and finite difference schemes \cite{FD_1,FD_2,FD_3,FD_4} are widely used due to their simplicity. The split-step Fourier method, in particular, leverages the fast Fourier transform to handle the linear and nonlinear parts of  the NLS equation separately, especially when dispersion and nonlinearity are well-balanced. Finite difference methods provide a straightforward grid-based approach for discretizing the NLS equation, though they can face accuracy issues for complex geometries. However, as the complexity of waveguide geometries has increased, along with the need for higher accuracy, finite element methods have gained prominence \cite{FEM_1,FEM_2,FEM_3,FEM_4}. Finite element methods offer greater flexibility in handling complex geometries, boundary conditions, and material properties, making them particularly suitable for inhomogeneous media and intricate waveguides.\\
A significant challenge in discretizing the NLS equation using finite element methods lies in formulating the problem as either a first-order or second-order system. While the second-order equation is stable and well-posed when discretized using least-squares finite element methods, it \mbox{often} leads to severe ill-conditioning, as the condition number grows at a rate of $\mathcal{O}(h^{-4})$ (see Theorem 9.11, \cite{LSFEM_a}), where $h$ is the mesh size. This ill-conditioning can significantly impair the accuracy of the approximate solution as meshes are refined. Conversely, the first-order system has a condition number that grows at a rate of $\mathcal{O}(h^{-2})$, which is better than $\mathcal{O}(h^{-4})$, but it leads to variational formulations that are inherently not well-posed, making it impossible to achieve stable and robust discretizations.\\
Thus, to address this issue, we examine the assumptions underlying the slowly varying envelope approximation \cite{NLOPT_book} used in deriving the NLS equation for pulse propagation in waveguides, and propose modifications to the existing NLS-based model. These modifications ensure that the resulting first-order system leads to a well-posed and stable variational formulation. Depending on the model parameters, the proposed model results in either a hyperbolic or an elliptic equation. In this article, we present a corresponding stability analysis of the proposed model for both hyperbolic and elliptic systems. Among the various finite element approaches, the use of the discontinuous Petrov-Galerkin (DPG) finite element method with optimal test functions is motivated by several key factors, such as mesh-independent stability, inheritance of discrete stability from the well-posedness of the continuous problem and the availability of a residual-based error estimator computed as part of the solution.\\
The article is organized as follows. Section~\ref{sec: laser_model} presents the detailed derivation of the laser propagation model. Section~\ref{sec: hyperbolic_case} and Section~\ref{sec: elliptic_judith_case} provide a detailed well-posedness analysis of the hyperbolic and elliptic system of equations, respectively. In Section~\ref{sec: numerical results}, we present numerical results for the proposed model for both hyperbolic and elliptic equations. Finally, we conclude with a short discussion in Section~\ref{sec: conclusions}.

\section{Pulsed laser propagation model} \label{sec: laser_model}
In this section, we provide a detailed description of the proposed model. The propagation of optical fields in fibers is governed by Maxwell's equations. The time dependent Maxwell's system of equations takes the form:
\begin{align}
\nabla \times \bm{E} &= -\frac{\partial \bm{B}}{\partial t} \label{mx_1}\\
\nabla \times \bm{H} &= \bm{J} + \frac{\partial \bm{D}}{\partial t} \label{mx_2}\\
\nabla \cdot \bm{D}  &= \rho_f \label{mx_3}\\
\nabla \cdot \bm{B}  &= 0\, \label{mx_4},
\end{align}
where $\bm{E}$ and $\bm{B}$ are the electric field and magnetic field vectors respectively, and $\bm{D}$ and $\bm{H}$ are the corresponding electric and magnetic flux densities. The current and free charge densities are represented by $\bm{J}$ and $\rho_f$, respectively. In silica fibers, we assume the absence of free charges, which leads to both the free charge and current densities being zero. $\bm{D}$ and $\bm{B}$ are related to  $\bm{E}$ and $\bm{H}$  by:
\begin{align}
\bm{D} = \varepsilon_0 \bm{E} + \bm{P}, \label{mx_5}\\
\bm{B} = \mu_0 \bm{H} + \bm{M}, \label{mx_6}
\end{align}
where $\bm{P}$ is the electric polarization vector, $\bm{M}$ is the magnetic polarization vector, $\varepsilon_0$ is the vacuum permittivity and $\mu_0$ is the vacuum permeability. For nonmagnetic media such as silica waveguides, magnetic polarization $\bm{M}$ can be neglected. Maxwell's equations can be utilized to obtain the wave equation for pulse propagation in optical waveguides. By taking the curl of~\cref{mx_1} and using~\cref{mx_2},~\cref{mx_5}, and ~\cref{mx_6}, $\bm{B}$ and $\bm{D}$ can be eliminated, yielding:
\begin{align}
\nabla \times \nabla \times \bm{E} = -\frac{1}{c^2} \frac{\partial^2 \bm{E}}{\partial t^2} - \mu_0 \frac{\partial^2 \bm{P}}{\partial t^2}, \label{curl_we}
\end{align}
where $c$ is the speed of light, and the relation $\mu_0 \varepsilon_0 = \frac{1}{c^2}$ is used. The polarization $\bm{P}$ can be expressed as a sum of linear and nonlinear contributions:
\begin{align}
\bm{P}(\bm{r},t) = \bm{P}_L(\bm{r},t) + \bm{P}_{NL}(\bm{r},t),
\end{align}
where $\bm{r} = {\left[x,y,z\right]}^T$. The linear contribution $\bm{P}_L$ and the nonlinear contribution $\bm{P}_{NL}$ are related to $\bm{E}$ by the following relations \cite{polarization_a,polarization_b,NLOPT_book_b}:
\begin{align}
     \bm{P}_L(\bm{r},t) &= \varepsilon_0 \int_{-\infty}^{t} \chi^{(1)}(t-t')\cdot \bm{E}(\bm{r},t')dt', \label{PL}\\
  \bm{P}_{NL}(\bm{r},t) &= \varepsilon_0 \int_{-\infty}^t \int_{-\infty}^t \int_{-\infty}^t \chi^{(3)}(t-t_1,t-t_2,t-t_3) \vdots \bm{E}  \otimes  \bm{E} \otimes \bm{E} \, dt_1dt_2dt_3, \label{PNL}
\end{align}
where $\chi^{(n)}$ represents the $n^{\text{th}}$ order susceptibility tensor, assuming the medium response is local. Under certain assumptions,~\cref{curl_we} can be further simplified. For the electric field vector $\bm{E}$, the following vector identity holds:
\begin{align}
\nabla \times \nabla \times \bm{E} = \nabla(\nabla \cdot \bm{E}) - \nabla^2\bm{E}. \label{VID}
\end{align}
In linear optics of isotropic, charge free materials, $\nabla \cdot \bm{E} = 0$ as $\nabla \cdot \bm{D} = 0$. However, in nonlinear optics, even for linear media, $\nabla \cdot \bm{E} \neq 0$, as a consequence of~\cref{PNL}.  The assumption $\nabla \cdot \bm{E} \approx 0$ is valid when the electric field takes the form of a transverse infinite wave \cite{NLOPT_book_b}. The contribution of $\nabla \cdot \bm{E}$ is shown to be small when the slowly-varying amplitude approximation is applicable (see Section 2.2,\cite{NLOPT_book_b}). With this simplification,~\cref{curl_we} reduces to:
\begin{align}
   -\nabla^2\bm{E} = -\frac{1}{c^2} \frac{\partial^2 \bm{E}}{\partial t^2} - \mu_0 \frac{\partial \bm{P}}{\partial t^2}. \label{lap_we}
\end{align}
\subsection{Nonlinear Pulse Propagation}
There are a few simplifying assumption that need to be made in order to solve~\cref{lap_we}. First, $\bm{P}_{NL}$ is treated as a small pertubation to $\bm{P}_L$, as nonlinear changes to the refractive index are less than $10^{-6}$. Second, the optical field is assumed to maintain its polarization along the length of the fiber. Third, the optical pulse is assumed to be quasi-monochromatic, with pulse being centered at $\omega_0$ and having a width of $\triangle \omega$ such that $\frac{\triangle \omega}{\omega_0} \ll 1$. In the slowly varying envelope approximation, the rapidly varying and slowly varying components of the electric field are separated, and the electric field is expressed explicitly as their product:
\begin{align}
    \bm{E}(\bm{r},t) = \frac{1}{2}\left[E(\bm{r},t)e^{-i\omega_0 t} + c\cdot c \right] \hat{x}, \label{ESV}
\end{align}
where $\hat{x}$ is the polarization unit vector, $E(\bm{r},t)$ is the slow varying component, the exponential represents the fast varying component, and $c \cdot c$ denotes the complex conjugate of the preceding term. Both $\bm{P}_L$ and $\bm{P}_{NL}$ can also be expressed in a similar way:
\begin{align}
        \bm{P}_L(\bm{r},t) = \frac{1}{2}\left[P_L(\bm{r},t)e^{-i\omega_0 t} + c\cdot c \right] \hat{x}, \label{PLSV}\\
        \bm{P}_{NL}(\bm{r},t) = \frac{1}{2}\left[P_{NL}(\bm{r},t)e^{-i\omega_0 t} + c\cdot c \right] \hat{x} \label{PNLSV},
\end{align}
where $P_L$ and $P_{NL}$ are the slow varying components. Another key assumption made while modeling pulse propagation is the instantaneous nonlinear response. This simplification reduces~\cref{PNL} to:
\begin{align}
    \bm{P}_{NL}(\bm{r},t) = \varepsilon_0 \chi^{(3)}\vdots \bm{E}(\bm{r},t)  \otimes  \bm{E}(\bm{r},t) \otimes \bm{E}(\bm{r},t). \label{PNLIR}
\end{align}
The assumption of an instantaneous nonlinear response amounts to neglecting the Raman effect. For silica fibers, the Raman response occurs over a timescale of $60$-$70$ps; hence, the assumption of an instantaneous nonlinear response is valid for pulses with a width greater than 1ps. Upon substituting~\cref{ESV} into~\cref{PNLIR}, $\bm{P}_{NL}$ is found to contain terms oscillating at $3\omega_0$ and $\omega_0$. The term corresponding to $3\omega_0$ can be neglected in optical fibers as it requires phase matching. By comparing the expression obtained by substituting~\cref{ESV} into~\cref{PNLIR} with~\cref{PNLSV}, we obtain:
\begin{align}
    P_{NL}(\bm{r},t) \approx \varepsilon_0\varepsilon_{NL} E(\bm{r},t), \label{PEF}
\end{align}
where $\varepsilon_{NL}$ represents the nonlinear contribution to the dielectric constant and is given by:
\begin{align}
    \varepsilon_{NL}(\bm{r},t) = \frac{3}{4} \chi^{3}_{xxxx} {\vert E(\bm{r},t) \vert}^2,\label{epnl}
\end{align}
where $\chi^{3}_{xxxx}$ represents the component of $\chi^3$ along the x-axis. While deriving the wave equation for pulse propagation, it is more convenient to use the slowly varying envelope approximation in the Fourier domain. Typically, this is not possible due to the presence of the nonlinear contribution to the dielectric constant. A widely used approach to address this issue is to treat $\varepsilon_{NL}$ as a constant during the derivation of the propagation equation \cite{DCNL_CONST}. This approach is justified due to the perturbative nature of  $P_{NL}$. Next, we define the Fourier transform of $E(\bm{r},t)$ as:
\begin{align}
    \Tilde{E}(\bm{r},\omega - \omega_0) = \int_{-\infty}^{\infty} E(\bm{r},t)e^{i(\omega - \omega_0)t}\, dt ,
\end{align}
and by substituting~\cref{ESV,PLSV,PNLSV,PNLIR,PEF,epnl} into~\cref{lap_we} and performing a Fourier transform, $\Tilde{E}(\bm{r},\omega - \omega_0)$ is found to satisfy:
\begin{align}
    \nabla^2 \Tilde{E} + \varepsilon(\omega) k_0^2 \Tilde{E} = 0, \label{Freq_eq}
\end{align}
where $k_0 = \frac{\omega}{c}$ and 
\begin{align}
    \varepsilon(\omega) = 1 + \Tilde{\chi}^{1}_{xx} + \varepsilon_{NL},
\end{align}
where $\Tilde{\chi}_{xx}^1$ denotes the Fourier transform of the component of $\chi^1$ aligned along the x-axis. The dielectric constant is expressed in terms of the modified refractive index \(\Tilde{n}\) and absorption coefficient \(\Tilde{\alpha}\). These modifications are given as follows \cite{NLOPT_book}:
\begin{align}
    \Tilde{n} = n + \frac{3}{8n} \mathcal{R}e\left(\chi^3_{xxxx}\right) {\vert E \vert}^2\quad \text{and} \quad \Tilde{\alpha} = \alpha + \frac{3\omega_0}{4nc} \mathcal{I}m \left(\chi^3_{xxxx}\right) {\vert E \vert}^2,
\end{align}
where $n$ and $\alpha$ are the linear refractive index and absorption coefficients, respectively.~\Cref{Freq_eq} can be solved by employing the method of separation of variables. Typically, $\Tilde{E}(\bm{r},\omega - \omega_0)$ is assumed to have the following form:
\begin{align}
    \Tilde{E}(\bm{r},\omega - \omega_0) = F(x,y)\Tilde{A}(z,\omega-\omega_0)e^{i\beta_0z}, \label{ansatz1}
\end{align}
where $\Tilde{A}$ is a slowly varying function of $z$, and $\beta_0$ is the wavenumber that will be determined later. On substituting~\cref{ansatz1} into~\cref{Freq_eq}, we obtain the following two equations for $F(x,y)$ and $\Tilde{A}$:
\begin{align}
    \frac{\partial^2 F}{\partial x^2} + \frac{\partial^2 F}{\partial y^2} + \left[ \varepsilon (\omega) k_0^2 - \Tilde{\beta}^2 \right]F = 0, \label{eq_mode}\\
    \frac{\partial^2 \Tilde{A}}{\partial z^2} + 2i \beta_0 \frac{\partial \Tilde{A}}{\partial z} + \left(\Tilde{\beta}^2 - \Tilde{\beta}_0^2\right) \Tilde{A} = 0, \label{eq_amp}
\end{align}
where  ${\Tilde{\beta}}^2$ is the separation constant. As a consequence of the slowly varying envelope approximation, $\frac{\partial^2 \Tilde{A}}{\partial z^2}$ is typically neglected while deriving the propagation equation.~\Cref{eq_mode} can be solved using first-order perturbation theory \cite{NLOPT_book_c}. In first-order perturbation theory, the nonlinear contribution to the dielectric constant doesn't affect the modal distribution, but the eigenvalue $\Tilde{\beta}$ becomes:
\begin{align}
    \Tilde{\beta}(\omega) = \beta(\omega) + \triangle \beta(\omega), \label{beta_pert}
\end{align}
where $\triangle \beta(\omega)$ represents the perturbation due to the nonlinearity in the dielectric constant. A closed-form expression for $\triangle \beta(\omega)$ can be found in \cite{NLOPT_book}, Section 2.3.1. Using~\cref{ansatz1} and~\cref{ESV}, the electric field $\bm{E}(\bm{r},t)$ can be written as 
\begin{align}
    \bm{E}(\bm{r},t) = \frac{1}{2}\left[F(x,y) A(z,t)e^{i(\beta_0z-\omega_0 t)} + c\cdot c \right] \hat{x},
\end{align}
where $A(z,t)$ is inverse Fourier transform of $\Tilde{A}(z,\omega-\omega_0)$. Typically, the exact form of $\beta(\omega)$ in~\cref{beta_pert} is unknown. Thus, in most propagation models, $\beta(\omega)$ is expanded around the carrier frequency $\omega_0$ using a Taylor series expansion as:
\begin{align}
    \beta(\omega) = \beta_0 + (\omega - \omega_0) \beta_1 + {(\omega - \omega_0)}^2 \beta_2 + \dots, \label{taylor}
\end{align}
where $\beta_i = \frac{\partial^i \beta}{\partial \omega^i}\big\rvert_{\omega = \omega_0}$ and $\beta_0 = \beta(\omega_0)$. A similar Taylor series expansion can be written for $\triangle \beta$ (Chapter 2,\cite{NLOPT_book}). For very short pulses with $\triangle \omega \ll \omega_0$, cubic and higher-order terms in the Taylor series expansion are negligible \cite{NLOPT_book}, and $\triangle \beta \approx \triangle \beta_0$ \cite{NLOPT_book}. By performing an inverse Fourier transform on~\cref{eq_amp} with these simplifications, we obtain the propagation equation for $A(z,t)$ in the time domain as:
\begin{align}
\frac{-i}{2\beta_0} \frac{\partial^2 A}{\partial z^2} + \frac{\partial A}{\partial z} + \beta_1 \frac{\partial A}{\partial t} + \frac{i \beta_2}{2} \frac{\partial^2 A}{\partial t^2} = i \triangle \beta_0 A. \label{TD_full}
\end{align}
In canonical models for pulse propagation, the term $\frac{-i}{2\beta_0} \frac{\partial^2 A}{\partial z^2}$ is neglected, as $A(z,t)$ is assumed to be a slowly varying function of $z$ as a consequence of the slowly varying envelope approximation. This simplifies the second-order equation to:
\begin{align}
    \frac{\partial A}{\partial z} + \beta_1 \frac{\partial A}{\partial t} + \frac{i \beta_2}{2} \frac{\partial^2 A}{\partial t^2} = i \triangle \beta_0 A. \label{CN_NLS}
\end{align}
\paragraph{Moving frame of reference:} In order to switch to a moving frame of reference propagating with the group velocity $v_g = \beta_1^{-1}$, a new coordinate system is defined as: 
\begin{align}
    \tau = t - \frac{z}{v_g},\, \text{and }  \xi = z. \label{mvfr}
\end{align}
Using the transformed coordinate system, we obtain the canonical nonlinear Schr\"{o}dinger equation for pulse wave propagation in the moving frame of reference:
\begin{equation}
    \frac{\partial A}{\partial \xi} + \frac{i \beta_2}{2}  \frac{\partial^2 A}{\partial \tau^2} = i \gamma {\vert A \vert}^2 A, \label{cn_NLS}
\end{equation}
where $ \triangle \beta_0 = \gamma {\vert A \vert}^2$ \cite{NLOPT_book}.
Let $\Omega_0 \subset \mathbb{R}$ be an open interval. The variable $\tau$ lies in $\Omega_0$, while the variable $\xi$ lies in the open interval $(0,Z)$ with $Z < \infty$. Next, let's consider the canonical Schr\"{o}dinger equation in~\cref{cn_NLS} with a general source term $f \in L^2(\Omega_0 \times (0,Z))$. The second-order NLS equation can be posed as the following first-order system:  
\begin{align}
    \sigma - \frac{\partial A}{\partial \tau} = g, & \quad \tau \in \Omega_0, 0 < \xi < Z,\\
    \frac{\partial A}{\partial \xi} + \frac{i \beta_2}{2}  \frac{\partial \sigma}{\partial \tau} = f,& \quad \tau \in \Omega_0, 0 < \xi < Z, \\
    A(\xi,\tau) = 0,&\quad \tau \in \partial \Omega_0, 0 < \xi < Z, \\
    A(0,\tau) = 0,&\quad \tau \in \Omega_0,
\end{align}
and formal equivalence between the first-order system and second-order equation when $g = 0$ can be claimed. However, the first-order system is not well-posed in $L^2(\Omega_0 \times (0,Z))$ (see \cite{Dem_NLS}). If the first-order system were well-posed, then there should exist constants $C_1,C_2 > 0$ such that:
\begin{equation*}
    {\Vert A \Vert}_{L^2(\Omega_0 \times (0,Z))} + {\Vert \sigma \Vert}_{L^2(\Omega_0 \times (0,Z))} \leq C_1{\Vert f \Vert}_{L^2(\Omega_0 \times (0,Z))} + C_2{\Vert g \Vert}_{L^2(\Omega_0 \times (0,Z))}. 
\end{equation*}
However, the second-order equation implies that:
\begin{equation*}
    {\big\Vert \frac{\partial A}{\partial \tau} \big\Vert}_{L^2(\Omega_0 \times (0,Z))}  = {\Vert g + \sigma \Vert}_{L^2(\Omega_0 \times (0,Z))}  \leq  C_1{\Vert f \Vert}_{L^2(\Omega_0 \times (0,Z))} + 2C_2{\Vert g \Vert}_{L^2(\Omega_0 \times (0,Z))},
\end{equation*}
for any solution $u$, which turns out to be false. A detailed exposition proving the fallacy of this claim can be found in \cite{Dem_NLS}. \\
Least-squares finite element methods can be used to produce stable discretizations of the second-order equation, but they will result in linear systems with a condition number that grows at a quartic rate \cite{LSFEM_a,LSFEM_b} with mesh refinements. This greatly hinders the use of mesh adaptivity to obtain high-fidelity solutions. In this article, we address this issue for wave propagation in optical waveguides by proposing certain modifications to the existing NLS equation-based model. These modifications allow the resulting second-order equation to be reformulated as a well-posed first-order system. Consequently, upon discretization using the DPG finite element method, the resulting linear system will exhibit a quadratic growth rate of the condition number. \\
In contrast to the canonical NLS-based propagation model, we propose that the term \(\frac{-i}{2\beta_0} \frac{\partial^2 A}{\partial z^2}\) in~\cref{TD_full} should not be neglected but retained as a perturbation to the existing NLS-based model. The modified propagation model, in the moving frame of reference, is expressed as:
\begin{align}
    \frac{-i}{2\beta_0} \left( \frac{\partial^2 A}{\partial \xi^2} + \beta_1^2 \frac{\partial^2 A}{\partial \tau^2} - 2 \beta_1 \frac{\partial^2 A}{\partial \xi \partial \tau}\right)+  \frac{\partial A}{\partial \xi} + \frac{i \beta_2}{2}  \frac{\partial^2 A}{\partial \tau^2} = i \gamma {\vert A \vert}^2 A. \label{mvref_eq1}
\end{align}
By multiplying~\cref{mvref_eq1} with $2\beta_0 i$ on both sides, we obtain:
\begin{align}
        \frac{\partial^2 A}{\partial \xi^2} - 2 \beta_1 \frac{\partial^2 A}{\partial \xi \partial \tau} -  \left(\beta_2 \beta_0 - \beta_1^2 \right)\frac{\partial^2 A}{\partial \tau^2} + 2\beta_0 i \frac{\partial A}{\partial \xi}= -2 \beta_0 \gamma {\vert A \vert}^2 A. \label{mvref_eq2}
\end{align}

\subsection{Nondimensionalization of the propagation model}
For silica waveguides and the fundamental mode, $\beta_0 = 4.44 \times 10^6 \, \text{m}^{-1}$, $\beta_1 = 104.01 \, \text{ps/m}$, and $\beta_2 = -0.447\, {\text{ps}}^2/\text{m}$. To nondimensionalize the equation, we select $10^{-10} \, \text{s}$ as the unit of time and $1 \, \text{m}$ as unit of space. This leads to $\mathcal{O}(\tau) \approx 10^{-10} \, \text{s}$ and $\mathcal{O}(\xi) \approx 1 \, \text{m} $, and we define  the variables, $\Bar{\tau}$ and $\bar{\xi}$ as $\Bar{\tau} = \tau/10^{-10} \, \text{s}$ and $\Bar{\xi} = \xi$. Upon substituting the nondimensionalized variables in~\cref{mvref_eq2}, we obtain:
\begin{align}
        \frac{\partial^2 A}{\partial \bar{\xi}^2} - 2 \bar{\beta}_1 \frac{\partial^2 A}{\partial \bar{\xi} \partial \bar{\tau}} -  \left(\bar{\beta}_2 \bar{\beta}_0 - \bar{\beta}_1^2 \right)\frac{\partial^2 A}{\partial \bar{\tau}^2} + 2\bar{\beta}_0 i \frac{\partial A}{\partial \bar{\xi}}= -2 \bar{\beta}_0 \gamma {\vert A \vert}^2 A, 
\end{align}
where $\bar{\beta}_0 \approx 10^{6}$, $\bar{\beta}_1 \approx 1$, and $\bar{\beta_2} \approx \pm 10^{-4}$. Therefore, we focus on solving the second-order PDE:
\be
u_{,\bar{\xi} \bar{\xi}} - 2 \bar{\beta}_1 u_{,\bar{\tau} \bar{\xi}} - \alpha u_{,\bar{\tau}  \bar{\tau}}  + 2 \bar{\beta}_0 i u_{,\bar{\xi}} = f(\bar{\tau}, \bar{\xi}), 
\label{eq:theproblem}
\ee
where $\alpha = \bar{\beta}_0 \bar{\beta}_2 - \bar{\beta}_1^2$ and $u$ represents $A(\bar{\tau},\bar{\xi})$. For positive $\bar{\beta}_2$, the equation is hyperbolic, while for the negative $\bar{\beta}_2$, the problem is elliptic. Next, we present the well-posedness analysis for the hyperbolic and the elliptic systems. In the analysis, we focus on the differential operator with source term $f(\bar{\tau}, \bar{\xi}) \in L^2(\Omega_0 \times (0,Z))$ being linear. The analysis for the nonlinear source requires further investigation, and will be presented in future work. In the analysis and numerical results presented from Section 3 onward, we will simplify the notation by dropping the bars from the non-dimensionalized variables. Additionally, for the sake of brevity, we will denote the domain $\Omega_0 \times (0,Z)$ as $D$, and will consider $\Omega_0 = [0,T]$. These changes will help maintain clarity without altering the analysis and conclusions.

\section{Well-Posedness Analysis: Hyperbolic Case} \label{sec: hyperbolic_case}
In this section, we present the well-posedness analysis for the hyperbolic system, i.e., when $\beta_2$ is positive, which is called \textit{normal} dispersion in the optics community. The analysis leverages the fact that a second-order hyperbolic equation can be transformed to an equivalent unique Friedrichs (symmetric) hyperbolic system of equations. If this system of equation meets certain conditions ensuring well-posedness. The first step involves reducing the second-order equation to a unique first-order Friedrichs system. We begin by postulating the first-order equations in the following form:
\begin{align}
    A u_{,\xi} + B v_{,\xi} + D u_{,\tau} + E v_{,\tau} & = 0,\label{Feq1} \\[8pt] 
    B u_{,\xi} + C v_{,\xi} + E u_{,\tau} + Fv_{,\tau}  &= 0\,.\label{Feq2}
\end{align}

By differentiating~\cref{Feq1} wrt. $\xi$ and~\cref{Feq2} wrt. $\tau$, and then summing them up, we obtain the following second-order equation:
\begin{equation}
    A u_{,\xi\xi} + (B+D) u_{,\tau \xi} + E u_{,\tau\tau} + B v_{,\xi\xi} + (C+E) v_{,\xi \tau} + F v_{,\tau\tau}  = 0 \, . \label{seq_a}
\end{equation}
On comparing the coefficients of differential operators from ~\cref{seq_a} with the coefficients of~\cref{eq:theproblem} to identify their correspondence, we obtain:
\begin{equation}
\left.
\begin{aligned}
B & = 0 \\
F & = 0\\
C+E & = 0 \\
A & = 1\\
E & = - \alpha \\
B+D & = - 2 \beta_1 
\end{aligned}
\right\}
\quad \Rightarrow \quad
\left\{ 
\begin{aligned}
A = 1 \\
B = 0 \\
C = \alpha \\
D = - 2 \beta_1 \\
E = - \alpha \\
F = 0.
\end{aligned}
\right. \label{seq_b}
\end{equation}
On substituting the values of coefficients from~\cref{seq_b} in~\cref{Feq1,Feq2}, we  arrive at the following first-order system:
\begin{equation}
\begin{pmatrix}
1 &  0 \\
0 & \alpha
\end{pmatrix}
\begin{pmatrix}
u_{,\xi}  \\ v_{,\xi} 
\end{pmatrix}
\, + \, 
\begin{pmatrix}
- 2 \beta_1 & - \alpha \\
- \alpha & 0 
\end{pmatrix}
\begin{pmatrix} u_{,\tau}  \\ v_{,\tau}  \end{pmatrix}
\, + \, 
\begin{pmatrix}
2 \beta_0 i & 0 \\
0 & 0 
\end{pmatrix}
\begin{pmatrix} u \\ v \end{pmatrix}
= 0 \, . \label{nonpositveeq}
\end{equation}
To demonstrate the well-possedness of the first-order hyperbolic system (\cref{nonpositveeq}), we verify whether the system of equations satisfies the assumptions (A1-A4) from Section 2.1 and assumptions (M1-M2) from Section 2.2 in \cite{Ern_Guermond_06}. The system of equations satisfies (A1-A3), but it fails to meet the postivity requirement in assumption A4 due to the presence of the zero-order term. To address this, we rescale the unknowns $(u,v)$ with an exponential factor:
\begin{equation}
    u := e^\xi u,\quad v := e^\xi v \, .
\end{equation}
This rescaling introduces the necessary zero-order term to satisfy the positivity condition. The rescaled equation is given by,
\begin{equation}
\underbrace{
\begin{pmatrix}
1 &  0 \\
0 & \alpha
\end{pmatrix}}_{=:\AB^\xi}
\begin{pmatrix}
u_{,\xi}  \\ v_{,\xi} 
\end{pmatrix}
\, + \, 
\underbrace{
\begin{pmatrix}
- 2 \beta_1 & - \alpha \\
- \alpha & 0 
\end{pmatrix}}_{:=\AB^{\tau}}
\begin{pmatrix} u_{,\tau}  \\ v_{,\tau}  \end{pmatrix}
\, + \, 
\begin{pmatrix}
1+2 \beta_0 i & 0 \\
0 & \alpha 
\end{pmatrix}
\begin{pmatrix} u \\ v \end{pmatrix}
= 0 \, . \label{positve_eq}
\end{equation}
We define a graph space $W$ as:
\begin{equation}
W = \{ \bm{w} \in {(L^2(D))}^2; \AB \bm{w} := \AB^\xi \bm{w}_{,\xi} + \AB^{\tau} \bm{w}_{,\tau} \, \in {(L^2(D))}^2\}, 
\end{equation}
equipped with the graph norm:
\begin{equation}
    {\Vert \bm{w} \Vert}_W^2 = {\Vert \AB \bm{w} \Vert}^2 + {\Vert \bm{w} \Vert}^2, 
\end{equation}
and associated inner product.
Now, let us consider a rectangular computational domain, as shown in~\Cref{fig:domain}, and the generalized eigenvalue problem associated with the hyperbolic system of equations:
\begin{figure}[h!]
	\centering
	\includegraphics[width=.5 \textwidth]{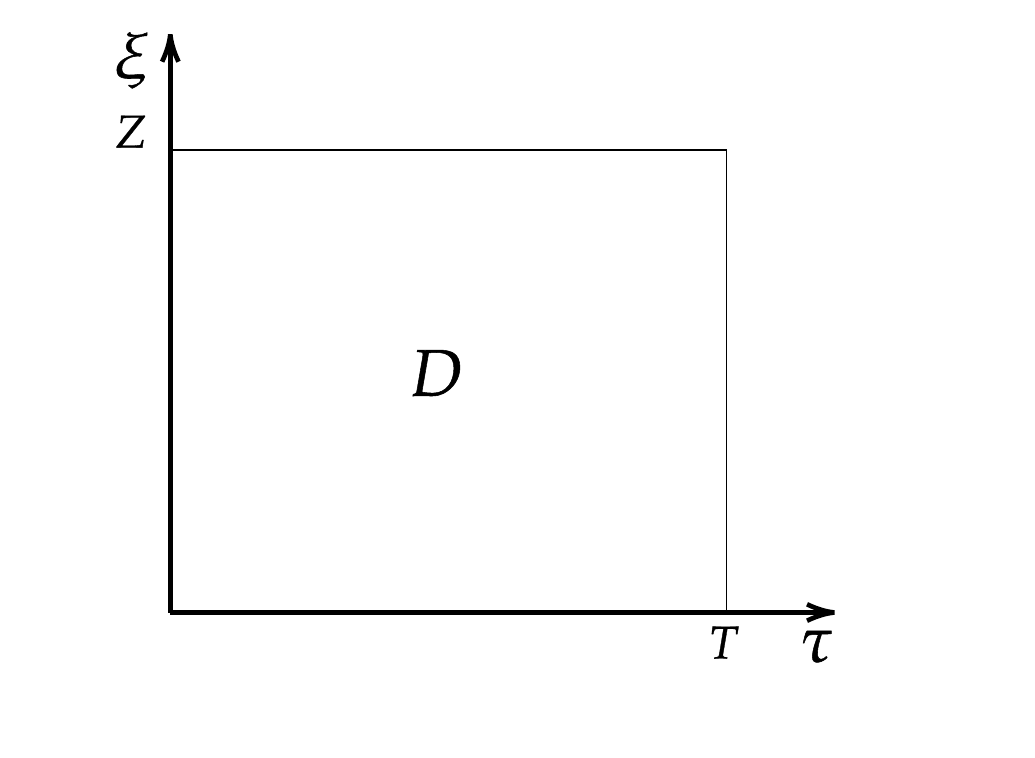}
	\caption{Domain of interest.}
	\label{fig:domain}
\end{figure}
\begin{equation}
\AB^\tau \bm{b} = \lambda \AB^\xi \bm{b} \, . \label{eigenvalueprob}    
\end{equation}
On performing spectral decomposition of~\cref{eigenvalueprob}, we obtain:
\begin{equation}
\begin{aligned}
    \lambda_1 = - \beta_1 - \sqrt{\beta_0 \beta_2}, \\
    \lambda_2 = - \beta_1 + \sqrt{\beta_0 \beta_2}, 
\end{aligned} \label{eigenvalues}
\end{equation}
as the eigenvalues, and 
\begin{equation}
    \begin{aligned}
        \bm{b}_1 = \frac{(- \lambda_1, 1)^T}{\sqrt{\lambda_1^2+ \alpha }}, \\
        \bm{b}_2 = \frac{(- \lambda_2, 1)^T}{\sqrt{\lambda_2^2 + \alpha}}, 
    \end{aligned}\label{eigenvectors}
\end{equation}
as the eigenvectors. On substituting the values of $\beta_0,\beta_1$ and $\beta_2$ in~\cref{eigenvalues} and~\cref{eigenvectors}, we obtain:
\begin{equation}
    \begin{aligned}
        \lambda_1 \approx -11,& \quad \bm{b}_1 \approx (0.740, 0.067)^T, \\
         \lambda_1 \approx 9,& \quad \bm{b}_2 \approx  (-0.669, 0.074)^T. 
    \end{aligned}
\end{equation}
The unit eigenvectors $\bm{b}_1,\bm{b}_2$ form an orthonormal basis with respect to the weighted inner product given by:
\begin{equation}
(\bm{U},\bm{V})_{\AB^\xi} :=  \bm{V}^\ast \AB^\xi \bm{U}  = u_1 \overline{v}_1 + \alpha u_2 \overline{v}_2 \, , \label{A_inner_product}
\end{equation}
where $\bm{U} = {[u_1,u_2]}^T  \, \in \, W$ and $\bm{V} = {[v_1,v_2]}^T \in W$. Next, one can represent the unknown $\bm{U}$ in terms of the eigenbasis:
\begin{equation}
\bm{U} = U_1 \bm{b}_1 + U_2 \bm{b}_2, \label{eigen_basis_rep}
\end{equation}
where $U_1$ and $U_2$ represents the corresponding components, and are given by:
\begin{equation}
\begin{aligned}
U_1 = (\bm{U}, \bm{b}_1)_{\AB^\xi} \approx \phantom{-}0.740 u + 6.727 v, \\[8pt]
U_2 = (\bm{U}, \bm{b}_2)_{\AB^\xi} \approx - 0.669 u + 7.432 v \, .
\end{aligned} \label{components}
\end{equation}

Using~\cref{A_inner_product} and~\cref{eigen_basis_rep}, we can express the transport contribution to the differential operator in~\cref{positve_eq} in terms of the eigenbasis and its components as follows: 
\begin{equation}
\begin{aligned}
\AB U = \AB^\tau \bm{U}_{,\tau} + \AB^\xi \bm{U}_{,\xi} & = U_{1,\tau} \AB^\tau  \bm{b}_1 + U_{2,\tau} \AB^\tau  \bm{b}_2 + U_{1,\xi} \AB^\xi  \bm{b}_1 + U_{2,\xi} \AB^\xi  \bm{b}_2 \\[8pt]
& = \underbrace{(\lambda_1 U_{1,\tau}  + U_{1,\xi} )}_{=: \ptl^{c_1} U_1} \AB^\xi  \bm{b}_1
+ \underbrace{(\lambda_2 U_{2,\tau}  + U_{2,\xi} )}_{=: \ptl^{c_2} U_2} A^\xi  \bm{b}_2, 
\end{aligned}
\end{equation}
where
\begin{equation}
 \ptl^{\bm{c}_1} U_1 = \underbrace{(\lambda_1, 1)}_{=: \bm{c}_1}  \cdot \nabla U_1 \qquad  \ptl^{\bm{c}_2} U_2 = \underbrace{(\lambda_2, 1)}_{=: \bm{c}_2}  \cdot \nabla U_2 \,.
\end{equation}
Hence, the $L^2$-integrability of $\AB \bm{U}$ is equivalent to the $L^2$-integrability of the directional derivatives $ \ptl^{\bm{c}_1} U_1 ,  \ptl^{\bm{c}_2} U_2$.
Additionally, on representing $\bm{V}$ in terms of its spectral components,
\begin{equation}
\bm{V} = V_1 \bm{b}_1 + V_2 \bm{b}_2,\qquad V_j = (\bm{V},\bm{b}_j)_{\AB^\xi},\, j=1,2\, ,
\end{equation}
we obtain:
\begin{equation}
(\AB \bm{U},\bm{V})_{L^2(D)} = (\lambda_1 U_{1,\tau} + U_{1,\xi}, V_1)_{L^2(D)} + (\lambda_2 U_{2,\tau} + U_{2,\xi}, V_2)_{L^2(D)} \, . \label{bilinear_form_a}
\end{equation}
By integrating~\cref{bilinear_form_a} by parts, we can define the boundary operator as:
\begin{equation}
\langle M \bm{U}, \bm{V} \rangle := \int_{\ptl D} (\bm{c}_1 \cdot \bm{n})_{+} U_1 \overline{V}_1 + (\bm{c}_2 \cdot \bm{n})_{+} U_2 \overline{V}_2 \, , \label{boundary_op}
\end{equation}
where $\bm{n}$ is the outward normal vector and $(\bm{c}_i \cdot \bm{n})_{+} \, \text{for} \, i=1,2$, represents positive contributions or outflow. The boundary operator $M$ satisfies the positivity assumption (see Section 2.2 in \cite{Ern_Guermond_06}):
\begin{equation}
    {\langle M \bm{U},\bm{U} \rangle}_{W,W'} \geq \,0 \quad\forall \, \bm{U} \in \, W,
\end{equation}
where $W'$ is dual space of $W$. This is a stronger assumption than the one made on the boundary operator in \cite{Ern_Guermond_06}. Therefore, the hyperbolic system satisfies all the assumptions made by the authors in \cite{Ern_Guermond_06}, namely assumptions A1–A4 and M1–M2, which are required for the well-posedness of the ultraweak variational formulation of the hyperbolic system (see Theorem 2.8 in \cite{Ern_Guermond_06}). Consequently, we define the stable variational formulation as:
\begin{equation}
    \text{find} \, \bm{U} \in W \, \text{such that} \, (\mathcal{A} \bm{U},\bm{V}) = (\bm{f},\bm{V}) \quad \forall \bm{V} \, \in \, W, 
\end{equation}
where $\bm{f} \in {(L^2(D))}^2$ and $\mathcal{A}$ is the differential operator from~\cref{positve_eq}.
To ensure that the trace of the convection operator is well-defined, a standard assumption is that the inflow and outflow boundaries must be separated by a positive distance. To satisfy this condition, we cut the corners of the rectangular domain along the characteristics.

\section{Well-Posedness Analysis for the Elliptic Case}\label{sec: elliptic_judith_case}
In this section, we examine the well-posedness of the ultraweak variational formulation of \eqref{eq:theproblem} for the elliptic case, i.e., when $\beta{_2}$ is negative, which is referred to as \textit{anomalous} dispersion in the optics community. For that, we first prove that the classical variational formulation for the second-order problem is well-posed, which allows us to conclude that the ultraweak formulation is also well-posed. 

\subsection{Second-order problem and classical variational formulation}
 Renaming $a := - \alpha > 0$, we rewrite equation \eqref{eq:theproblem} in divergence form with the following boundary conditions
\begin{equation}\label{eq:divform}
\displaystyle{ \left\{
\begin{split}
-\div(\bm{A}\nabla u)-2\beta_0iu_{,\xi}&=f, & \,\mbox{in } D,\\
u&=u_{0}, &\, \mbox{ on }\Gamma_1,\\
(\bm{A}\nabla u)\cdot \bm{n}&=t, & \, \mbox{ on }\Gamma_2.\\
\end{split}
\right.} 
\end{equation}
Here, $n$ denotes the unit outward vector on $\Gamma=\Gamma_1\cup\Gamma_2$, the diffusion tensor is
\begin{align}
    \bm{A} = \left(
\begin{array}{cc}
a & - \beta_1 \\
- \beta_1 & 1
\end{array}
\right),
\end{align}
$\Gamma_1$ is the part of the boundary consisting of three edges where $\xi=0$ or $\tau=\left\{0,T\right\}$ and $\Gamma_2$ is the top boundary, i.e., where $\xi=Z$. 

For the well-posedness analysis, it is sufficient to consider homogeneous boundary conditions $u_0=t=0$. Introducing the following energy space
$$
V := \{ v \in H^1(D) \, : v =0 \text{ on } \Gamma_1 \}\, ,
$$
we now consider the classical variational formulation
\be
\left\{
\begin{array}{ll}
u \in V \\[8pt]
\ds \int_D (\bm{A}\nabla u)\cdot\nabla \bar{v}- 2 \beta_0 i u_{,\xi} \bar{v} = \int_D f \bar{v},\;\;\forall v \in V \, .
\end{array}
\right.
\label{eq:VF_2nd_order_elliptic}
\ee
We need to prove that the operator induced by the sesquilinear form in \eqref{eq:VF_2nd_order_elliptic} is bounded below. 

To this end, first note that the leading term in \eqref{eq:VF_2nd_order_elliptic} is elliptic because the diffusion tensor $\bm{A}$ is positive-definite. By standard reasoning involving the Poincaré inequality (see \cite{FE_Math_book_SIAM}, Section 2.3.1.), it is easy to see then that the first term in \eqref{eq:VF_2nd_order_elliptic} is coercive.

The first-order term in \eqref{eq:VF_2nd_order_elliptic} represents the compact perturbation of a coercive term, and therefore, the Fredholm Alternative applies (see e.g. \cite{FAbook}, Section 5.21). The latter concludes that in order to prove that the operator is bounded below, it is sufficient to demonstrate uniqueness.
Considering the homogeneous case $f=0$ and testing with $v=u$ in \eqref{eq:VF_2nd_order_elliptic}, we obtain
\begin{equation}\label{eq:VF_injectivity}
\int_D (\bm{A}\nabla u)\cdot\nabla \bar{u}- 2 \beta_0 i u_{,\xi} \bar{u} =\underbrace{\int_D (a u_{,\tau} - \beta_1 u_{,\xi}) \bar{u}_{,\tau} + (- \beta_1 u_{,\tau} + u_{,\xi} ) \bar{u}_{,\xi}}_{\text{real number}}   - 2 \beta_0 i \int_D  u_{,\xi} \bar{u}= 0.
\end{equation}
Note that the first integral on the right-hand side of \eqref{eq:VF_injectivity} is a real number, so in order for the whole integral to vanish, the imaginary part in the second term must be zero. Thus, the second term is re-written as 
$$
\begin{array}{ll}
\ds i \int_D u_{,\xi}\bar{u}= i \int_D (\Re u_{, \xi} + i \Im u_{,\xi} )(\Re u - i \Im u)
& \ds = i  \int_D  \half [ ( \Re u )^2 + (\Im u )^2 ]_{, \xi}  + \int_D \Re u_{, \xi} \Im u - \Im u _{,\xi} \Re u \\[12pt]
& \ds = \underbrace{i \int_0^T   \half \lvert u(\tau,Z) \rvert^2  \, d \tau}_{\text{imaginary number}} +  \underbrace{\int_D \Re u_{, \xi} \Im u - \Im u _{,\xi} \Re u }_{\text{real number}}\, .
\end{array}
$$
Setting the imaginary part to zero, one obtains that $u = 0$ on the top boundary, i.e, at $\xi=Z$. At the same time, we have from the second boundary condition in \eqref{eq:divform} that $-\beta_1u_{,\tau}+u_{,\xi}=0$ also at $\xi=Z$. These two conditions imply that $u_{,\xi}=0$ along the top boundary as well. Therefore, we can extend the solution $u(\tau,\xi)$ with zero to a larger domain, say $(0,T)\times(0,2Z)$. But the
solution to the elliptic problem with constant coefficient is analytic and, by analytic continuation argument, u must be zero in the whole domain.

In conclusion, the uniqueness of the solution follows and, therefore, the second-order formulation is well-posed. The well-posedness in $V \subset H^1(D)$ implies that we can replace $f$ with a general  $l \in V^\prime$, being $V'$ the dual space of $V$, and claim the existence
of a positive constant $C>0$ such that
$$
\Vert u \Vert_{H^1(D)} \leq  C \, \Vert l \Vert_{V'} \, .
$$ 

\subsection{First-order system and ultraweak variational formulation}
Introducing an extra variable, the flux $\sigma$, the first-order system reads 

\begin{equation}\label{eq:elliptic_1st_order_system}
\displaystyle{ \left\{
\begin{split}
\bm{\sigma}-\bm{A}\nabla u &=\bm{g}, & \mbox{in } D,\\
-\div\bm{\sigma}-2\beta_0iu_{,\xi}&=f, & \mbox{in } D,\\
u&=u_{0}, & \, \mbox{ on }\Gamma_1,\\
\bm{\sigma}\cdot \bm{n}&=t, & \,\mbox{ on }\Gamma_2.\\
\end{split}
\right.} 
\end{equation}
with $\bm{g}\in (L^2(D))^2$ and $f\in L^2(D)$. Note that for the original problem \eqref{eq:divform} we have that $\bm{g}=[0,0]^T$. The operator governing the first-order system is
\begin{equation}\label{eq:strongoperator}
\mathcal{A}(\bm{\sigma},u):=(\bm{\sigma}-\bm{A}\nabla u,-\div\bm{\sigma}-2\beta_0iu_{,\xi}),    
\end{equation}
with domain 
$$\mathcal{D}(\mathcal{A}):=\left\{(\bm{\sigma},u)\in (L^2(D))^2\times L^2(D)\,:\,\mathcal{A}(\bm{\sigma},u)\in (L^2(D))^2\times L^2(D),\, u=0\mbox{ on }\Gamma_1,\,\bm{\sigma}\cdot \bm{n}=0\mbox{ on }\Gamma_2\right\}.$$
Following the arguments from  \cite{Demkowicz_15}, we will use the Closed Range Theorem to study the well-posedeness of the first-order system \eqref{eq:elliptic_1st_order_system}. 
For that, we need to prove boundedness below of operator $\mathcal{A}$ defined in \eqref{eq:strongoperator}, determine its adjoint $\mathcal{A}^*$ and its null space. In order to show the boundedness below, we need to prove that $u$ and $\sigma$ are controlled by the $L^2$-norms of the righ-hand-side.

Again, we consider the case with homogeneous boundary condition, i.e., $u_0 = t = 0$. Let $v \in V$ be an arbitrary test function, we multiply the first equation with $-\nabla \bar{v}$ and the second equation with $\bar{v}$. Integrating over $D$,
relaxing the second equation 
and summing it all up,  we obtain the original second-order equation \eqref{eq:VF_2nd_order_elliptic} for $u$
with the right-hand side equal to 
$$
l(v) := - (\bm{g}, \nabla v) + (f,v) \,,
$$
where $(\cdot,\cdot)$ denotes either the usual inner product in $L^2(D)$ or $(L^2(D))^2$. The well-posedness of the second-order problem implies then that
$$
\Vert u \Vert ,\,\Vert \nabla u \Vert \lesssim \Vert l \Vert_{V^\prime} \leq  (\Vert \bm{g} \Vert + \Vert f \Vert )
$$
where $\Vert \cdot \Vert$ denotes the norms in $L^2(D)$ and $(L^2(D))^2$.  Moving the terms with derivatives of $u$ to the right-hand side in \eqref{eq:elliptic_1st_order_system}, and using the triangle
inequality, we obtain bounds:
$$
\Vert \bm{\sigma} \Vert ,\, \Vert \div \bm{\sigma} \Vert \lesssim (\Vert \bm{g} \Vert + \Vert f \Vert) \,.
$$

Using direct integration-by-parts, we conclude that the adjoint of operator $\mathcal{A}$ is 
\begin{equation}\label{eq:adjointoperator}
\mathcal{A}^*(\bm{\tau},v):=(\bm{\tau}+\nabla v,\div(\bm{A}\bm{\tau})-2\beta_0iv_{,\xi}),   
\end{equation}
and incorporating the boundary conditions in the definition of the domain, we have
\begin{equation}\nonumber
\begin{split}
\mathcal{D}(\mathcal{A}^*):=\left\{(\bm{\tau},v)\in (L^2(D))^2\times L^2(D)\,:\,\mathcal{A}^*(\bm{\tau},v)\in (L^2(D))^2\times L^2(D),\,v=0\mbox{ on }\Gamma_1,\right.\\
\left.-(\bm{A}\bm{\tau})\cdot \bm{n} -2\beta_0 iv=0\,\mbox{ on }\Gamma_2\right\}.
\end{split}
\end{equation}
If we eliminate $\bm{\tau}$ to obtain a single equation for $v$, we obtain system \eqref{eq:divform} with the same essential boundary conditions on $\Gamma_1$ but the natural boundary condition of $\Gamma_2$ reads
$$(\bm{A}\nabla{v})\cdot \bm{n}-2\beta_0iv=0.$$
Repeating {\em exactly the same argument} as for the original variational formulation~(\ref{eq:VF_2nd_order_elliptic}),
the uniqueness for $v$ is demonstrated, from which it follows that the adjoint of the operator \eqref{eq:adjointoperator}
is injective as well.

We conclude then that the strong formulation for the first-order system is well-posed, and the Closed Range argument (see \cite{Demkowicz_15} for details) implies that
the ultraweak formulation of~(\ref{eq:elliptic_1st_order_system}) is likewise well-posed.
\section{DPG discretization} \label{sec:DPG_disc}

The core essence of the DPG methodology is the automatic generation of a stable discretization for a given well-posed variational formulation and an approximate trial space. The method attains stability by approximating an optimal discrete test space, tailored to the variational problem and the trial space, ensuring that the supremum in the discrete inf-sup condition in the continuous setting is automatically attained over the discrete test space. In the ideal DPG method, the optimal test space is obtained by inverting the Riesz operator over a discontinuous test space. However, in practice, the inversion of Riesz operator is not feasible due to its infinite dimensional nature. Thus, the inverse of the Riesz operator is approximated by inverting a Gram matrix induced by the test norm over an enriched finite dimensional test space. The discontinuous nature of test space allows local element-wise inversion of the Gram matrix.

\subsection{Hyperbolic System}

\begin{figure}[htbp]
	\centering
	\includegraphics[width=.4 \textwidth]{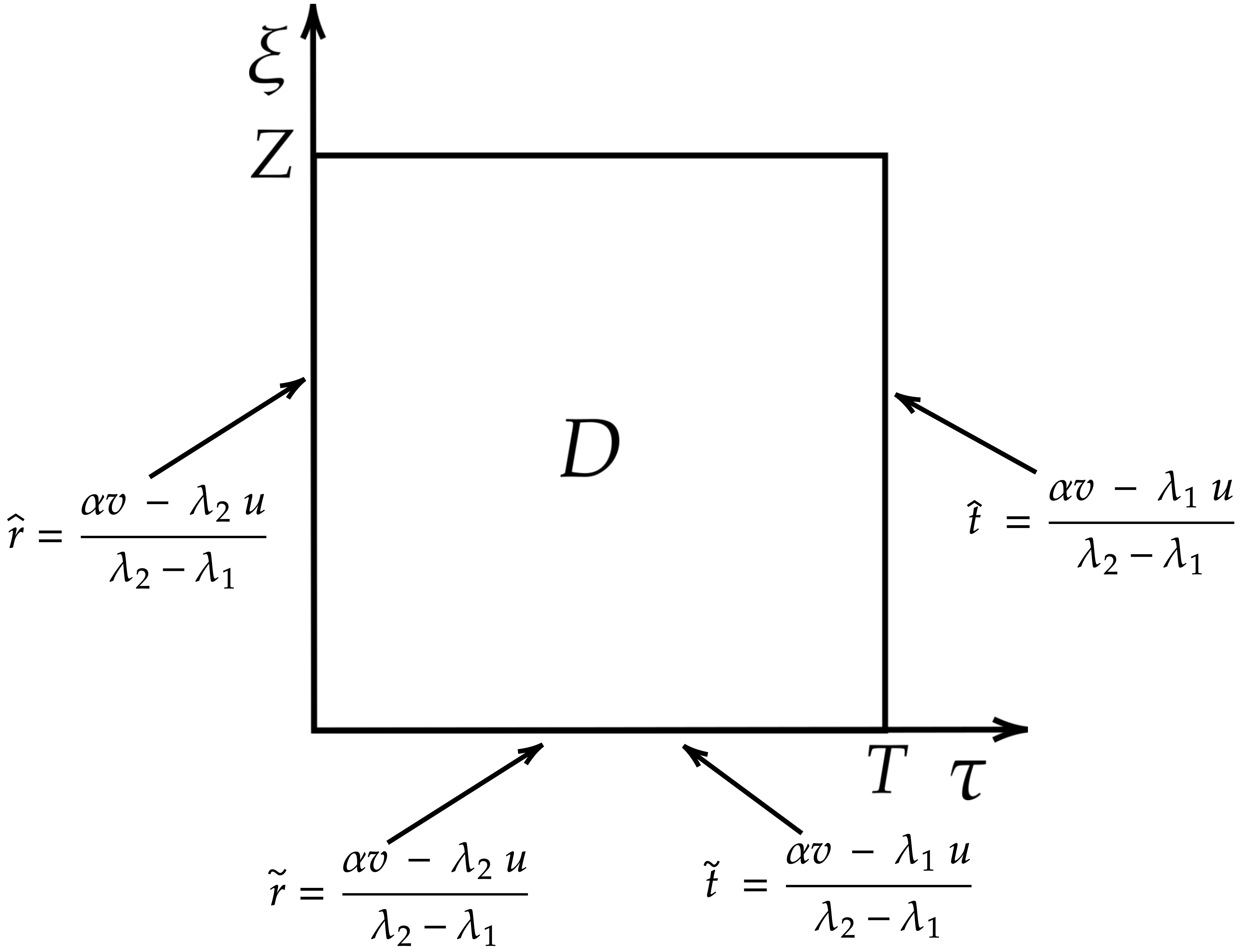}
	\caption{Boundary conditions for the hyperbolic system of equations.}
	\label{fig:hyp_bc}
\end{figure}

\begin{figure}[htbp]
	\centering
	\includegraphics[width=.4 \textwidth]{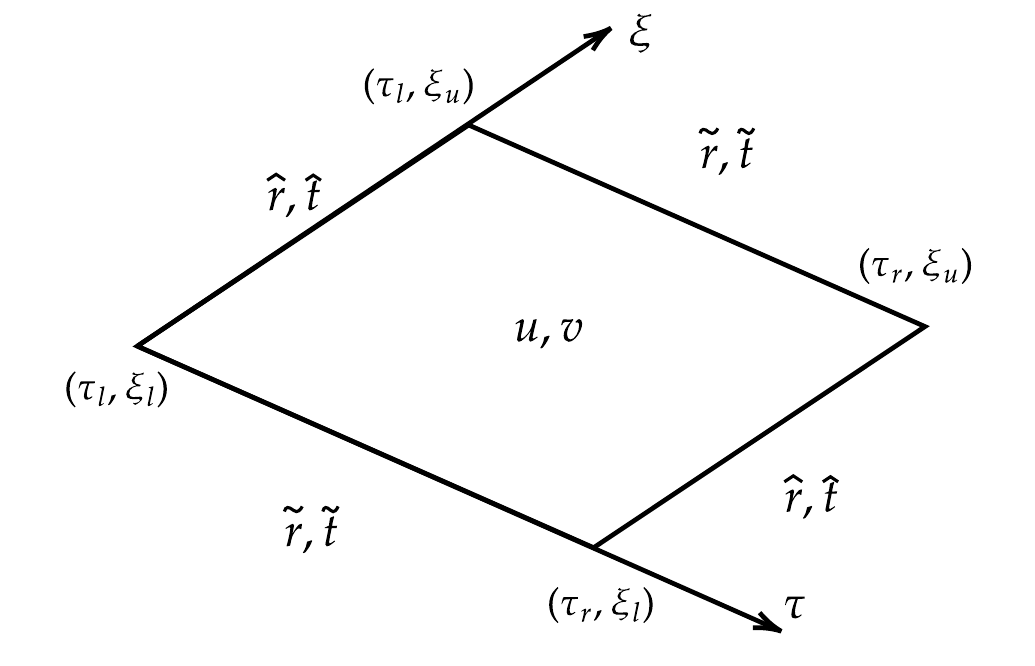}
	\caption{A tensor product element $K$ for the hyperbolic problem.}
	\label{fig:hyp_Tp_elem}
\end{figure}

Let's consider the system of hyperbolic equation from~\cref{positve_eq} with a source vector $\bm{f} ={[f_2(\xi,\tau), \\f_1(\xi,\tau)]}^T \in {(L^2(D))}^2$ and the boundary operator stated in~\cref{boundary_op}. The boundary conditions are obtained by using the negative definite part of the boundary integral:
\begin{equation}
     {\langle M \bm{U},\bm{V} \rangle}_{W,W'} = \int_{\ptl D} (\bm{c}_1 \cdot \bm{n})_{-} U_1 \overline{V}_1 + (\bm{c}_2 \cdot \bm{n})_{-} U_2 \overline{V}_2, 
\end{equation}
where $(\bm{c}_i \cdot \bm{n})_{-} \, \text{for} \, i=1,2$, represents the negative contributions along $\bm{c}_1$ and $\bm{c}_2$. For a square domain $D={[0,1]}^2$, only the domain boundaries defined by $\xi = 0$, $\tau = 0$, and $\tau = 1$ will have negative contributions from $(\bm{c}_i \cdot \bm{n})U_i  \, \, \text{for} \, \, i\in {1,2}$. These negative contributions are the quantities that are needed to be prescribed as the inflow boundary conditions (see~\Cref{fig:hyp_bc}).
Using~\cref{positve_eq} and the boundary conditions from~\Cref{fig:hyp_bc}, we define the hyperbolic problem as:
\begin{align*}
\begin{pmatrix}
1&0\\0&\alpha
\end{pmatrix} \begin{pmatrix}
u_{,\xi} \\ v_{,\xi}
\end{pmatrix} + \begin{pmatrix}
-2 \beta_1 & -\alpha \\ -\alpha & 0
\end{pmatrix}  \begin{pmatrix}
u_{,\tau} \\ v_{,\tau}
\end{pmatrix} + \begin{pmatrix}
2\beta_0 i&0\\0&0
\end{pmatrix}\begin{pmatrix}
u \\ v
\end{pmatrix} &= \begin{pmatrix}
f_1(\xi,\tau) \\ f_2(\xi,\tau)
\end{pmatrix} & \quad \forall (\tau,\xi) \in  \Omega_0 \times I, \\
\frac{\alpha v - \lambda_1 u}{\lambda_2 - \lambda_1} &= g_1(\tau) &\quad \forall (\tau,\xi) \in  \Omega_0 \times \{ 0 \},\\
\frac{\alpha v - \lambda_2 u}{\lambda_2 - \lambda_1} &= g_2(\tau) &\quad \forall (\tau,\xi) \in  \Omega_0 \times \{0\}, \\
\frac{\alpha v - \lambda_2 u}{\lambda_2 - \lambda_1} &= g_3(\xi) &\quad \forall (\tau,\xi) \in  \{0\} \times I,\\
\frac{\alpha v - \lambda_1 u}{\lambda_2 - \lambda_1} &= g_4(\xi) &\quad \forall (\tau,\xi) \in  \{T\} \times I,
\end{align*}
where $I = [0,Z]$. We define tensor product energy spaces for the hyperbolic system as: 
\begin{equation}
\begin{aligned}
H^1(\Omega_0) \otimes L^2(I) &:=\{ q: \Omega_0 \times I \rightarrow \mathbb{R} : \left({\Vert q(\xi,\tau) \Vert}^2_{\Omega_0 \times I} + {\Vert  q_{,\tau}(\xi,\tau)\Vert}^2_{\Omega_0 \times I} \right) < \infty \},   \\
L^2(\Omega_0) \otimes H^1(I) &:= \{ r: \Omega_0 \times I \rightarrow \mathbb{R} : \left({\Vert r(\xi,\tau) \Vert}^2_{\Omega_0 \times I} +  {\Vert r_{,\xi}(\xi,\tau)\Vert}^2_{\Omega_0 \times I} \right) < \infty \}, \\
L^2(\Omega_0) \otimes L^2(I) &:= \{ s: \Omega_0 \times I \rightarrow \mathbb{R} : {\Vert s(\xi,\tau) \Vert}^2_{\Omega_0 \times I}  < \infty \}. 
\end{aligned}
\end{equation}
In the DPG methodology, discontinuous energy spaces are used for the test functions. Hence, we define the broken equivalent of  $H^1(\Omega_0 \times I)$ for the tensor product finite element mesh $(D_h)$:
\begin{align*}
H^1(D_h) &:= \{ u:\Omega_0 \times I \rightarrow \mathbb{R}: u\vert_K \in H^1(\Omega_{0,K} \times I_K) \,\, \forall \,\, K \, \in \, D_h\},
\end{align*}
where $K \in D_h$ represents an element of the tensor product finite element mesh. We present an illustration of one such tensor product element in~\Cref{fig:hyp_Tp_elem}.
A major implication of using broken test spaces is the introduction of additional trace variables on the mesh skeleton. The trace spaces are defined as:
\begin{align*}
\hat{\Lambda} &:= \{ \hat{\lambda} : \exists \, \lambda \in H^1(\Omega_0) \otimes L^2(I) \, \text{such that} \, \hat{\lambda} = \gamma^K_{\tau}(\lambda\vert_K)  \, on \, \partial K_{\xi} \quad \forall \, K \, \in \, D_h  \},  \\
\tilde{\Lambda} &:=\{ \tilde{\lambda} : \exists \, \lambda \in L^2(\Omega_0) \otimes H^1(I) \, \text{such that} \, \tilde{\lambda} = \gamma^K_{\xi}(\lambda\vert_K)  \, on \, \partial K_{\tau} \quad \forall \, K \, \in \, D_h  \},
\end{align*}
where $\gamma^K_{\tau}$ and $\gamma^K_{\xi}$ are continuous trace operators along $\tau$ and $\xi$ directions respectively, $\partial K_{\xi}$ and $\partial K_\tau$ represents the boundary of the tensor product element $K$ along $\xi$ and $\tau$ axis.
\paragraph{Ultraweak formulation:} Let $(U,\hat{U})$ be the approximate trial space, $V$ the test space, and $V'$ the dual space of V. The ultraweak DPG formulation of the hyperbolic system can be stated as: Given $l \in V'$, find $\mathfrak{u} \in U$ and $\hat{\mathfrak{u}} \in \hat{U}$ satisfying:
\begin{equation}
    b(\mathfrak{u},\mathfrak{v}) + \hat{b}(\hat{\mathfrak{u}},\mathfrak{v}) = l(\mathfrak{v}) \quad \forall \mathfrak{v} \in V, \label{eq:blf_1}
\end{equation}
where
\begin{equation}
    \begin{aligned}
        \mathfrak{u} &= (u,v) \in U := L^2(D) \times L^2(D), \\
        \hat{\mathfrak{u}} &= (\hat{r},\hat{t},\tilde{r},\tilde{t}) \in \hat{U} := \hat{\Lambda} \times \hat{\Lambda} \times \tilde{\Lambda} \times  \tilde{\Lambda}, \\
        \tilde{r} &= g_1(\tau) \quad \forall (\tau,\xi) \in  \Omega_0 \times \{ 0 \},\\
        \tilde{t} &= g_2(\tau) \quad \forall (\tau,\xi) \in  \Omega_0 \times \{0\}, \\
        \hat{r}   &= g_3(\xi)  \quad \forall (\tau,\xi) \in  \{0\} \times I, \\
        \hat{t}   &= g_4(\xi)  \quad \forall (\tau,\xi) \in  \{T\} \times I, \\
        \mathfrak{v} &=(\delta u, \delta v) \in V: =  H^1(D_h) \times H^1(D_h) \\
        b(\mathfrak{u},\mathfrak{v}) &=  \sum_{K \in D_h}-{(u,{\delta u}_{,\xi})}_K + 2\beta_1{(u ,{\delta u}_{,\tau})}_K + \alpha {(v ,{\delta u}_{,\tau})}_K + 2\beta_0 i{(u ,\delta u)}_K - \alpha {(v ,{\delta v}_{,\xi})}_K + \alpha {(u ,{\delta v}_{,\tau})}_K,   \\
        \hat{b}(\hat{\mathfrak{u}},\mathfrak{v}) &= \sum_{K \in D_h}{\int_{\tau_l}^{\tau_r} (\tilde{t} - \tilde{r}) \delta u \, d\tau} \bigg\rvert_{\xi = \xi_l}^{\xi = \xi_u} -2\beta_1  {\int_{\xi_l}^{\xi_u} (\hat{t} - \hat{r}) \delta u \, d\xi} \bigg\rvert_{\tau = \tau_l}^{\tau = \tau_r}-{\int_{\xi_l}^{\xi_u} (\lambda_2\hat{t} - \lambda_1\hat{r}) \delta u \, d\xi} \bigg\rvert_{\tau = \tau_l}^{\tau = \tau_r},\\
        &+{\int_{\tau_l}^{\tau_r} (\lambda_2\tilde{t} - \lambda_1\tilde{r}) \delta v \, d\tau} \bigg\rvert_{\xi = \xi_l}^{\xi = \xi_u} -\alpha  {\int_{\xi_l}^{\xi_u} (\hat{t} - \hat{r}) \delta v \, d\xi} \bigg\rvert_{\tau = \tau_l}^{\tau = \tau_r}, \\
        l(\mathfrak{v}) &= \sum_{K \in D_h}(f_1,\delta u)_K + (f_2,\delta v)_K. 
    \end{aligned} \label{eq:wf_hyp}
\end{equation}
The broken test spaces are equipped with the adjoint graph norm \cite{Dem_adj_grp,Dem_chan}:
\begin{equation}
    {\Vert \mathfrak{v} \Vert}^2 := {\Vert \mathcal{A}_h^{\star} \mathfrak{v} \Vert}^2 + \alpha { \Vert \mathfrak{v} \Vert}^2,
\end{equation}
where $\alpha > 0$ is a scaling constant, and $\mathcal{A}_h^{\star}$ is (formal) adjoint operator computed element-wise. The adjoint (formal) operator is given by:
\begin{equation}
    \mathcal{A}_h^{\star}(\delta u, \delta v) = (-{\delta u}_{,\xi} + 2\beta_1 {\delta u}_{,\tau} + 2\beta_0 i {\delta u} + \alpha {\delta v }_{,\tau}, -\alpha {\delta v}_{,\xi} + \alpha {\delta u}_{,\tau})_{D_h}.
\end{equation}
In this paper, all experiments are done using $\alpha = 1$.
\subsection{Elliptic System}
Next, let's consider the elliptic system of equations given in~\cref{eq:elliptic_1st_order_system}. From~\cref{eq:elliptic_1st_order_system}, it can be inferred that the leading order term in the first equation ($\sigma - A^{-1}\sigma = g$) is of $\mathcal{O}(1)$, while the leading order term in the second equation $(-\div \bm{\sigma} - 2\beta_0 iu_{,\xi }= f)$ is of $\mathcal{O}(10^6)$ . This large difference in scale in the equations of the first-order elliptic system leads to extremely ill-conditioned linear system obtained upon discretization (see~\cref{err_scaling_comp}). Hence, we rescale $\bm{\sigma}$ in the first-order system with a scaling factor $c$, such that each equation in the first-order system attains a comparable scale. The rescaled first-order elliptic system \footnote{In~\cref{eq:rescaled_firstorder}, the negative sign in front of the differential operator: $\div\bm{\sigma}+2c^{-1}\beta_0iu_{,\xi}$ has been absorbed in the source term $f$.} reads
\begin{equation}\label{eq:rescaled_firstorder}
\displaystyle{ \left\{
\begin{split}
c\bm{\sigma}-\bm{A}\nabla u&=\bm{g},&\mbox{in } D,\\
\div\bm{\sigma}+2c^{-1}\beta_0iu_{,\xi}&=c^{-1}f,&\mbox{in } D,\\
u&=u_{0},&\mbox{ on }\Gamma_1,\\
\bm{\sigma}\cdot \bm{n}&=t,&\mbox{ on }\Gamma_2.\\
\end{split}
\right.}
\end{equation}
Next, we define the energy spaces required for the elliptic problem as:
\begin{align}
L^2(\Omega_0 \times I) &= \left\{ u : \Omega_0 \times I \rightarrow \mathbb{R}: {\Vert u \Vert} < \infty \right\}, \nonumber \\
H^1(\Omega_0 \times I) &= \left\{ v : \Omega_0 \times I \rightarrow \mathbb{R}:  v \, \in \, L^2(\Omega_0 \times I), \nabla v \, \in \, {\left(L^2(\Omega_0 \times I)\right)}^2   \right\}, \\
\bm{H}(\div,\Omega_0 \times I) 
&= \left\{  \bm{w}: \Omega_0 \times I \rightarrow \mathbb{R}^2: \bm{w} \, \in \, (L^2(\Omega_0 \times I))^2, \nabla \cdot \bm{w} \, \in \, L^2(\Omega_0 \times I)	\right\}. \nonumber
\end{align}
The broken equivalents of $H^1(\Omega_0 \times I)$ and $\bm{H}(\div,\Omega_0 \times I)$ are defined as:
\begin{align}
    H^1(D_h) &:= \left\{ v : \Omega \rightarrow \mathbb{R}: v\big|_{K} \in H^1(K) \quad \forall \, K \, \in \, D_h  \right\}, \\
\bm{H}(\div,D_h) &:=  \left\{ \bm{w} : \Omega \rightarrow \mathbb{R}^3: \bm{w}\big|_{K} \in \bm{H}(\div,K) \quad \forall \, K \, \in \, D_h 	\right\},
\end{align}
where $K$ represents an element in the finite element mesh $D_h$.  The trace spaces are defined as:
\begin{align}
\begin{split}
H^{1/2}(\Gamma_h) &:= \left\{ \hat{u} : \exists \, u \, \in \, H^1(\Omega_0 \times I) \, \text{such that} \, \hat{u} = \gamma^K(u\big|_{K}) \, \text{on} \, \partial K \quad \forall \, K \in D_h  \right\}, \\
H^{-1/2}(\Gamma_h) &:= \left\{ \hat{\sigma}_n : \exists \, \bm{\sigma} \, \in \bm{H}(\div,\Omega_0 \times I)\, \text{such that} \, \hat{\sigma}_n = \gamma_n^K(\bm{\sigma} \big|_{K}) \, \text{on} \, \partial K \quad \forall \, K \in D_h  \right\},
\end{split}
\end{align}
where  $\gamma^K$ and $\gamma_n^K$ represent continuous and normal trace operators, respectively \cite{lec_LD}.
\paragraph{Ultraweak formulation:} Let $(U,\hat{U})$ be the trial space, $V$ the test space, and $V'$ the dual space of V. The ultraweak DPG formulation of the elliptic system can be stated as: Given $l \in V'$, find $(\mathfrak{u} \in U$ and $\hat{\mathfrak{u}} \in \hat{U}$ satisfying:
\begin{equation}
    b(\mathfrak{u},\mathfrak{v}) + \hat{b}(\hat{\mathfrak{u}},\mathfrak{v}) = l(\mathfrak{v}) \quad \forall \mathfrak{v} \in V, \label{eq:blf_2}
\end{equation}
where
\begin{equation}
    \begin{aligned}
        \mathfrak{u} &= (u,\bm{\sigma}) \, \in \, U :=\ L^2(\Omega_0 \times I) \times (L^2(\Omega_0 \times I))^2 , \\
\hat{\mathfrak{u}} &= (\hat{u},\hat{\sigma}_n) \, \in \, \hat U :=\ H^{1/2}(\Gamma_h) \times H^{-1/2}(\Gamma_h): \hat{u} = u_0 \, \text{on} \, \Gamma_1, \hat{\sigma}_n = \sigma_0 \, \text{on} \, \Gamma_2, \\
\mathfrak{v} &= (v,\bm{\tau}) \, \in \, V :=\ H^1(D_h) \times \bm{H}(\div,D_h), \\
b(\mathfrak{u},\mathfrak{v}) &= -(\bm{\sigma},\nabla v)_{D_h} - 2c^{-1}\beta_0 i(u,v_{,\xi}) + (u,\div \bm{\tau}) + c(\bm{A}^{-1} \bm{\sigma},\bm{\tau})\\
\hat{b}(\hat{\mathfrak{u}},\mathfrak{v}) &=2c^{-1}\beta_0 i \langle \bm{n}_{\xi}\hat{u},v \rangle\ + \langle \hat{\sigma}_n,v \rangle - \langle \hat{u},\bm{\tau} \cdot \bm{n} \rangle\ \\
l(\mathfrak{v}) &= (\bm{g},\tau)_{D_h} + (f,v)_{D_h} .
    \end{aligned}  \label{eq:wf_ellp}
\end{equation}
The broken adjoint (formal) operator for elliptic system is given by:
\begin{equation}
    \mathcal{A}_h^{\star}(v, \bm{\tau}) = (2c^{-1}\beta_0i v_{,\xi} + \div \bm{\tau}, -\nabla v + c \bm{A}^{-1}\bm{\tau}).
\end{equation}
\paragraph{Implementation of DPG method:} A typical viewpoint --- in terms of implementation of the DPG methodology --- is to interpret the problems given by \cref{eq:blf_1,eq:blf_2} as a mixed method \cite{DPG_ency,Dahmen_dpg_1}. Let $V_h \subset V$ denote the enriched test space, and let:
\begin{equation}
    R_{V,K} : V(K) \rightarrow (V(K))'
\end{equation}
denote the element-wise Riesz operator corresponding to the inner product associated with the test space. By using the fact that $R_{V,K}$ is an isometry, we can identify the residual defined in the broken test space $\psi_k :=  R_{V,K}^{-1}(l(\cdot) - b(\mathfrak{u}_h,\cdot) - \hat{b}(\hat{\mathfrak{u},\cdot})$ as a new unkown, and solve the following mixed problem:
\begin{equation}
\left\{
	\arraycolsep=2pt
	\begin{array}{lll}
		\text{find }\psi_h \in V_h,\, \mathfrak{u}_h \in U_h,\, \hat{\mathfrak{u}}_h \in U_h , \\
		(\psi_h, \mathfrak v)_V - b(\mathfrak u_h, \mathfrak v) - \hat b(\hat{\mathfrak u}_h, \mathfrak v) 
		&= l(\mathfrak v) & \quad \mathfrak v \in V_h, \\
		b(\delta \mathfrak u_h, \psi_h) 
		&= 0 & \quad \delta \mathfrak u_h \in U_h, \\
		\hat b(\delta \hat{\mathfrak u}_h, \psi_h) 
		&= 0 & \quad \delta \hat {\mathfrak u}_h \in \hat U_h.
	\end{array}
 \right.
\end{equation}
Next, we briefly discuss the algebraic structure of the resulting linear system and the built-in error estimator. Let the basis functions of $V_h,U_h$ and $\hat{U}_h$ be denoted by $\varphi_i,\psi_i$ and $\hat{\psi}_i$ respectively. From~\cref{eq:wf_hyp} for the hyperbolic system or~\cref{eq:wf_ellp} for the elliptic system, we can construct the following matrices for each element $K \in D_h$:
\begin{equation}
\begin{split}
{\mathrm{G}}_{K,lj} &= {(\varphi_l,\varphi_j)}_{V} \, , \\[8pt]
\mathrm{B}_{K,ij} &= b_K(\varphi_i,\psi_j) \, , \\[8pt]
\mathrm{\hat{B}}_{K,ij} &= \hat{b}_K(\varphi_i,\hat{\psi}_j) \, , \\[8pt]
\mathrm{l}_{K,i} &= l_K(\varphi_i),
\end{split} \label{components_of_linear_system}
\end{equation}
where ${\mathrm{G}}_{K,lj}$ represents the element Gram matrix corresponding to the test inner product and approximates the Riesz operator when discretized over $V_h(K)$. $\mathrm{B}_{K,ij} $ represents the element stiffness matrix corresponding to the $L^2$ variables,  $\mathrm{\hat{B}}_{K,ij}$ represents the element stiffness matrix corresponding to the trace variables,  and $\mathrm{l}_{K,i}$ is the element load vector. Using~\cref{components_of_linear_system}, we obtain a symmetric positive definite linear system for element K,
\begin{equation}
{\begin{bmatrix}
\mathrm{G}_K & \mathrm{B}_{K} & \mathrm{\hat{B}}_{K} \\
\mathrm{B}_{K}^T & \mathrm 0 & \mathrm 0 \\
\mathrm{\hat B}_{K}^T & \mathrm 0 & \mathrm 0
\end{bmatrix}}
\begin{bmatrix} 
\mathrm{\Psi}_K \\
\mathrm{u}_K \\
\hat{\mathrm{u}}_K
\end{bmatrix}
= 
\begin{bmatrix}
\mathrm{l}_K \\
\mathrm 0 \\
\mathrm 0 \\
\end{bmatrix} \, ,
\end{equation}
where  $[\mathrm{u}_K, \hat{\mathrm{u}}_K]^T$ and $[\mathrm{\Psi}_K]$  are the solution DOFs and DOFs associated with the residual respectively. In practical implementations, the residual dofs are statically condensed, resulting in a linear system of the form:
\begin{equation}
{\begin{bmatrix}
\mathrm{B}_{K} & \mathrm{\hat{B}}_{K}
\end{bmatrix}}^T
\mathrm{G}_K^{-1} 
\begin{bmatrix}
\mathrm{B}_{K} & \mathrm{\hat{B}}_{K}
\end{bmatrix} 
\begin{bmatrix} 
\mathrm{u}_K & \hat{\mathrm{u}}_K
\end{bmatrix}^T
= 
{\begin{bmatrix}
\mathrm{B}_{K} & \mathrm{\hat{B}}_{K}
\end{bmatrix}}^T 
\mathrm{G}_K^{-1}
\mathrm{l}_K \, .
\end{equation}
Once the solution DOFs are computed by solving the globally assembled system, the element-wise DPG residual is computed as:
\begin{equation}
\mathrm{\Psi}_K =
\mathrm{G}^{-1}_K 
(\mathrm{l}_K - \mathrm{B}_K \mathrm{u}_K - \mathrm{\hat{B}}_K \mathrm{\hat u}_K ) \, .
\end{equation}
An in-depth exposition of the algebraic structure of the linear system generated by the ultraweak DPG finite element formulations can be found in \cite{Dem_part_2,Dem_polyDPG}.
\section{Numerical Results} \label{sec: numerical results}

In this section, we present numerical results motivated by our interest in pulse propagation in optical waveguides. Optical solitons play a crucial role in pulse propagation, especially in optical fiber communications. A soliton is stable, self-reinforcing wave packet that maintains its shape while travelling at constant velocity, due to the balance between the dispersion and the nonlinear effects in the medium. 
These solitons allow transmission of data over long distances without distortion as they resists pulse broadening. 

\begin{figure}[htbp]
	\centering
    \begin{subfigure}[t]{0.45\textwidth}

        \includegraphics[width=\textwidth]
        {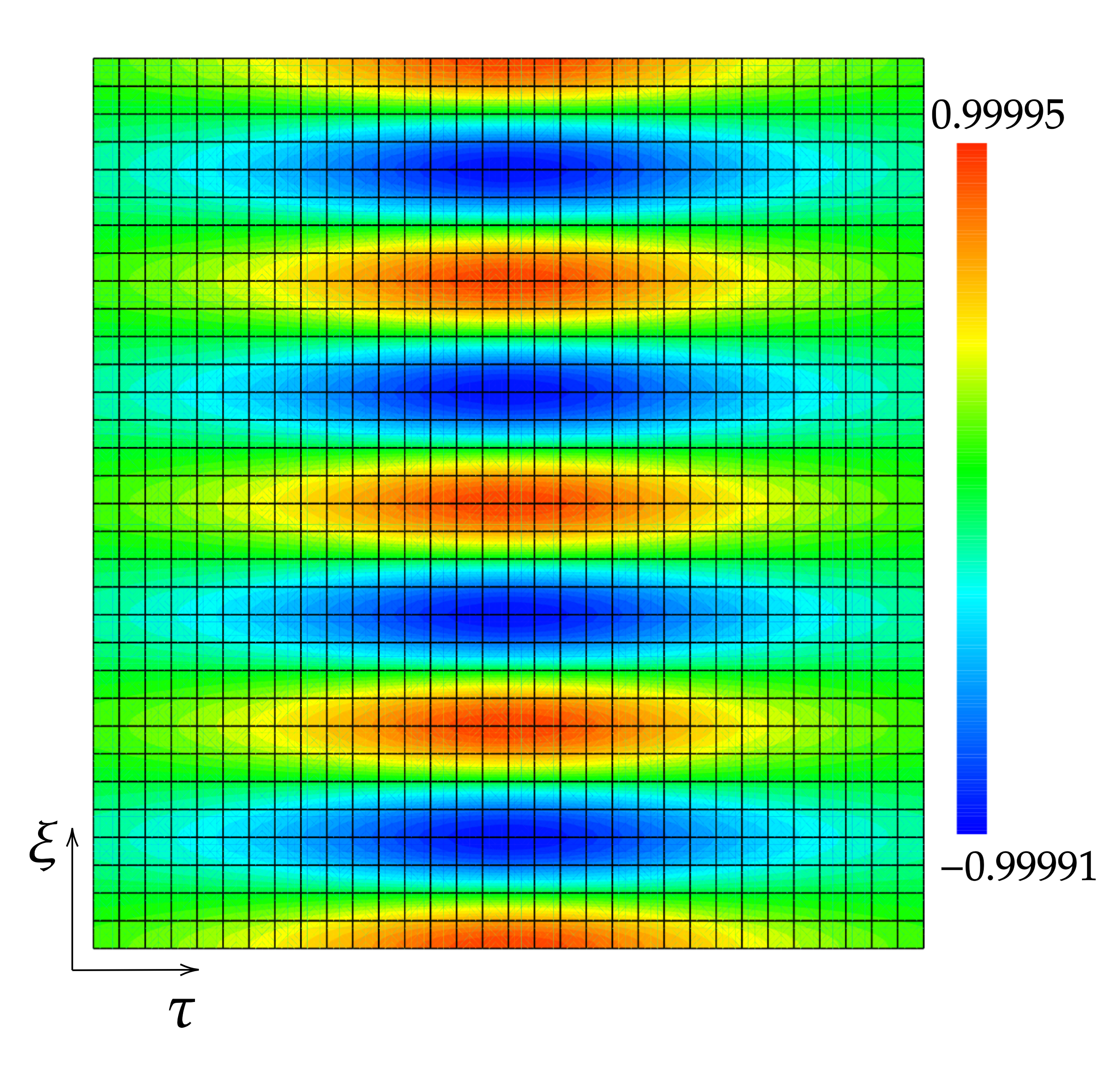}
        		\caption{\small $\Re(u(\tau,\xi))$}
    \end{subfigure}
	\begin{subfigure}[t]{0.45\textwidth}

		\includegraphics[width=\textwidth]
		{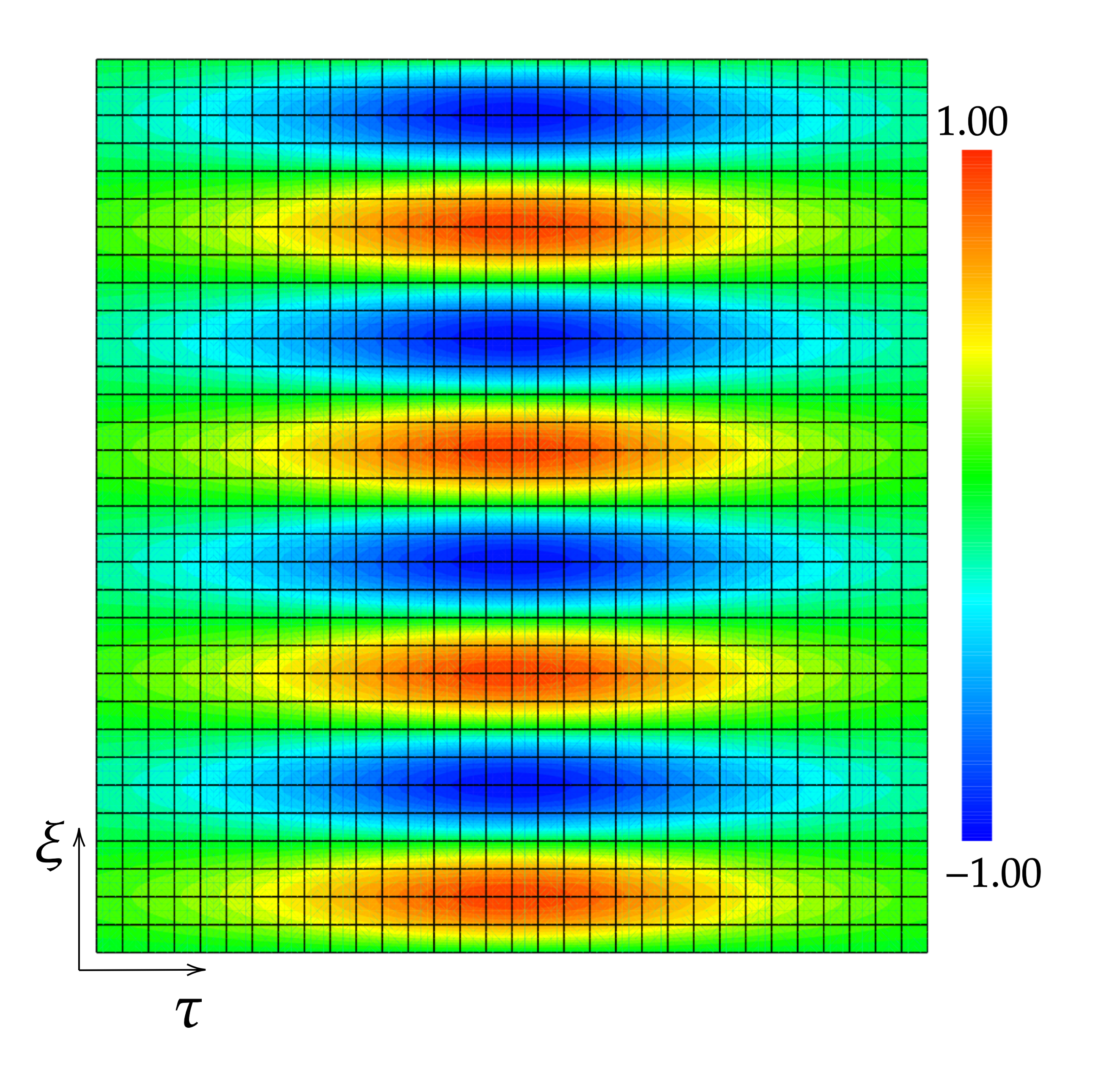}
  		\caption{\small $\Im(u(\tau,\xi))$}
    \end{subfigure}
    \caption{First-order soliton: Contour plots of the (a) real and (b) imaginary part of the solution $u$ for the hyperbolic system.}\label{solFOS}
\end{figure}

\begin{figure}[htbp]
	\centering
    \begin{subfigure}[t]{0.45\textwidth}

        \includegraphics[width=\textwidth]
        {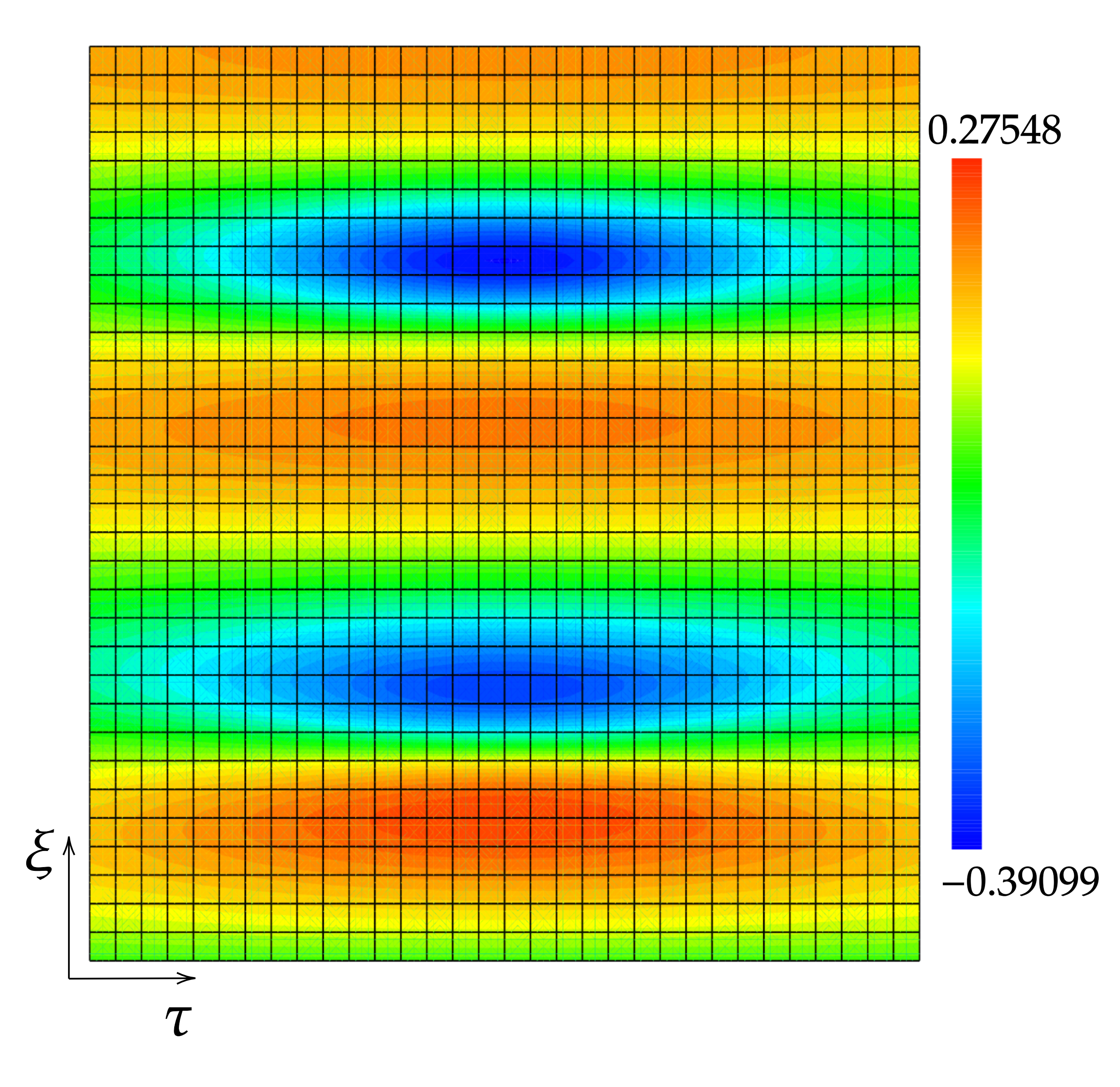}
        		\caption{\small $\Re(u(\tau,\xi))$}
    \end{subfigure}
	\begin{subfigure}[t]{0.45\textwidth}

		\includegraphics[width=\textwidth]
		{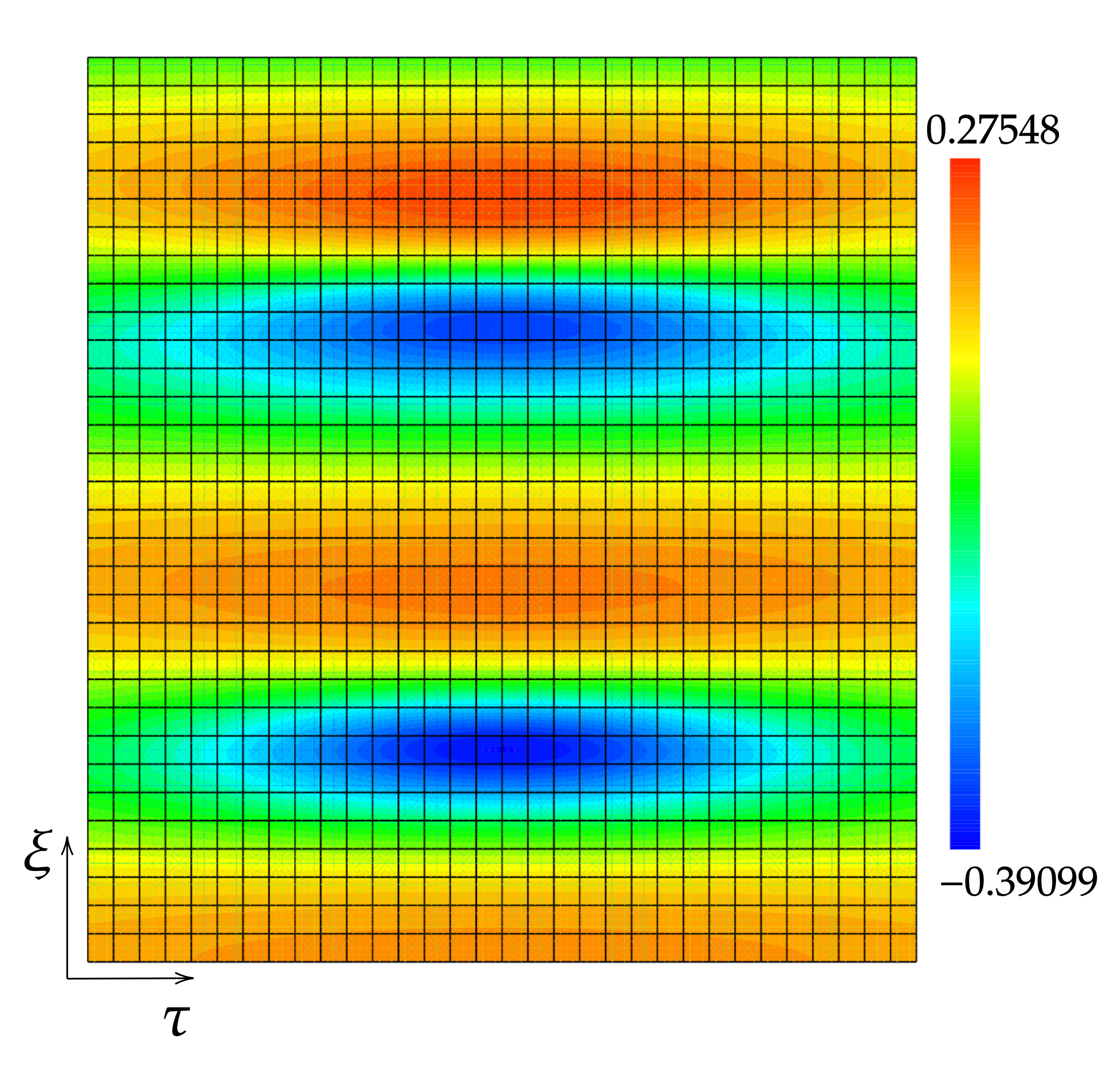}
  		\caption{\small $\Im(u(\tau,\xi))$}
    \end{subfigure}
    \caption{Second-order soliton: Contour plots of the (a) real and (b) imaginary part of the solution $u$ for the hyperbolic system.}\label{solSOS}
\end{figure}

We exemplify the success of our methodology by comparing against exact soliton solutions. 
A first-order soliton solution in an optical waveguide can be formulated as
\begin{align}
    u(\xi,\tau) = \mathrm{sech}(a_0(\tau-0.5))e^{i \omega \xi},
\end{align}
where $\omega = 8\pi$ and $a_0 = 5.0$. Likewise, a second-order soliton can be expressed as 
\begin{align}
    u(\tau,\xi) = \frac{4(\cosh(3(\tau-0.5)) + 3\cosh(\tau-0.5)e^{4i\xi})e^{\frac{i \omega \xi}{2}}}{\cosh(4(\tau-0.5))+4\cosh(2(\tau-0.5))+3\cos(4 \omega \xi)},
\end{align}
where $\omega = \pi$. 
For reference, the solutions of first-order and second-order solitons are depicted in~\Cref{solFOS,solSOS}, respectively. 

Next, we establish the appropriate the source terms ($f(\xi,\tau)$) and boundary conditions in the first-order elliptic and hyperbolic systems in preparation for the numerical solves. 
The first-order scaled elliptic system is 
\begin{equation}\label{eq:rescaled_firstorder}
\displaystyle{
\begin{split}
c\bm{\sigma}-\bm{A}\nabla u&=\bm{g},&\mbox{in } D,\\
\div\bm{\sigma}+2c^{-1}\beta_0iu_{,\xi}&=c^{-1}f,&\mbox{in } D,\\
u&=u_{0},&\mbox{ on }\Gamma_1,\\
\bm{\sigma}\cdot \bm{n}&=t,&\mbox{ on }\Gamma_2,\\
\end{split}
}
\end{equation}
where
\begin{align}
    \bm{A} = \left(
\begin{array}{cc}
\alpha & - \beta_1 \\
- \beta_1 & 1
\end{array}
\right)
\quad
\text{and} \quad
\bm{A}^{-1} = \frac{1}{\alpha - \beta_1^2}
\left(
\begin{array}{cc}
1 &  \beta_1 \\
 \beta_1 & \alpha
\end{array}
\right)
\, .
\end{align}
And, the first-order hyperbolic system is 
\begin{align*}
\begin{pmatrix}
1&0\\0&\alpha
\end{pmatrix} \begin{pmatrix}
u_{,\xi} \\ v_{,\xi}
\end{pmatrix} + \begin{pmatrix}
-2 \beta_1 & -\alpha \\ -\alpha & 0
\end{pmatrix}  \begin{pmatrix}
u_{,\tau} \\ v_{,\tau}
\end{pmatrix} + \begin{pmatrix}
2\beta_0 i&0\\0&0
\end{pmatrix}\begin{pmatrix}
u \\ v
\end{pmatrix} &= \begin{pmatrix}
f_1(\xi,\tau) \\ f_2(\xi,\tau)
\end{pmatrix} \quad  \text{in} \quad D, \\
\frac{\alpha v - \lambda_1 u}{\lambda_2 - \lambda_1} &= g_1(\tau) \quad \text{on} \quad  [0,1] \times \{ 0 \},\\
\frac{\alpha v - \lambda_2 u}{\lambda_2 - \lambda_1} &= g_2(\tau) \quad \text{on} \quad    [0,1] \times \{0\}, \\
\frac{\alpha v - \lambda_2 u}{\lambda_2 - \lambda_1} &= g_3(\xi) \quad \text{on} \quad    \{0\} \times (0,1],\\
\frac{\alpha v - \lambda_1 u}{\lambda_2 - \lambda_1} &= g_4(\xi) \quad \text{on} \quad    \{T\} \times (0,1].
\end{align*}
Additionally, we prescribe $v = \cos(\tau) + i \sin(\xi)$ to compute the boundary conditions and source vector for the hyperbolic problem. 

\begin{figure}[htbp]
\centering
\resizebox{0.4\textwidth}{!}{
\begin{tikzpicture}
		\begin{loglogaxis}[xmin=1,xmax=150, ymin=2e-3,ymax=200,xlabel=\large{$\sqrt{\text{\# Elements}}$},ylabel=\large{$\text{Relative } L^2 \text{ error}$},grid=major,legend style={at={(1,1)},anchor=north east,font=\small,rounded corners=2pt}]
		\addplot[color = blue,mark=square*] table[x= sqrt_n, y=rel_L2_error, col sep = comma] {results/conv_p1_SOS2_hyp_new_para.dat};
  		\addplot[dashed,color = blue,forget plot] table[x= sqrt_n, y=extslp, col sep = comma] {results/conv_p1_SOS2_hyp_new_para.dat};
		\addplot [color = red,mark=square*] table[x= sqrt_n, y=rel_L2_error, col sep = comma] {results/conv_p2_SOS2_hyp_new_para.dat};
    		\addplot[dashed,color = red,forget plot] table[x= sqrt_n, y=extslp, col sep = comma] {results/conv_p2_SOS2_hyp_new_para.dat};
		\addplot [color = black,mark=square*] table[x= sqrt_n, y=rel_L2_error, col sep = comma] {results/conv_p3_SOS2_hyp_new_para.dat};
    		\addplot[dashed,color = black,forget plot] table[x= sqrt_n, y=extslp, col sep = comma] {results/conv_p3_SOS2_hyp_new_para.dat};
		\legend{$p = 2$,$p = 3$,$p = 4$}
		\end{loglogaxis}
	\end{tikzpicture}
}
\hspace{0.05\textwidth}
\resizebox{0.4\textwidth}{!}{
\begin{tikzpicture}
		\begin{loglogaxis}[xmin=1,xmax=150, ymin=1e-4,ymax=100,xlabel=\large{$\sqrt{\text{\# Elements}}$},ylabel=\large{$\text{Residual}$},grid=major,legend style={at={(1,1)},anchor=north east,font=\small,rounded corners=2pt} ]
		\addplot[color = blue,mark=square*] table[x= sqrt_n, y=res, col sep = comma] {results/conv_p1_SOS2_hyp_new_para.dat};
		\addplot [color = red,mark=square*] table[x= sqrt_n, y=res, col sep = comma] {results/conv_p2_SOS2_hyp_new_para.dat};
		\addplot [color = black,mark=square*] table[x= sqrt_n, y=res, col sep = comma] {results/conv_p3_SOS2_hyp_new_para.dat};
		\legend{$p = 2$,$p = 3$,$p = 4$}
		\end{loglogaxis}
\end{tikzpicture}
}
\caption{Second-order soliton: convergence plots for the (a) relative $L^2$ error in $\Re(u(\tau,\xi))$ and $\Re(v(\tau,\xi))$, and (b) DPG residual for the hyperbolic system.} \label{errConvSOShyp}
\end{figure}
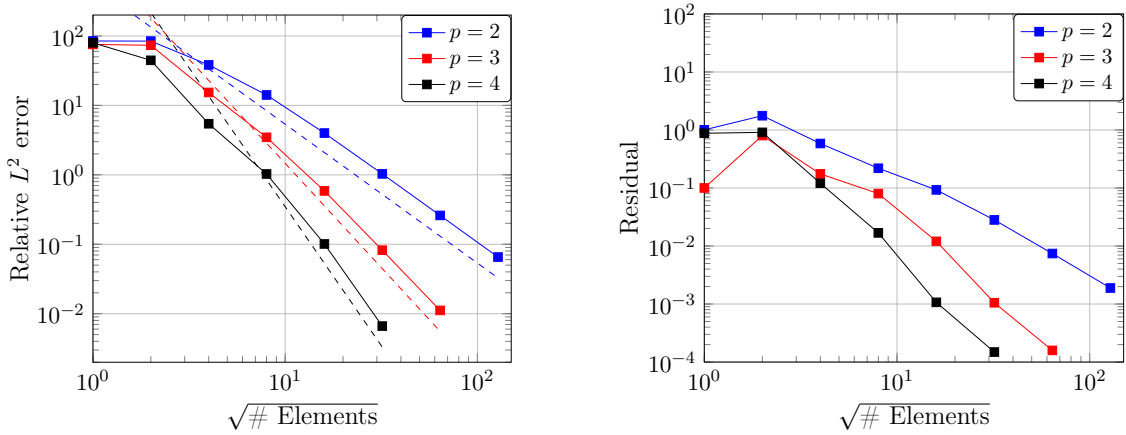

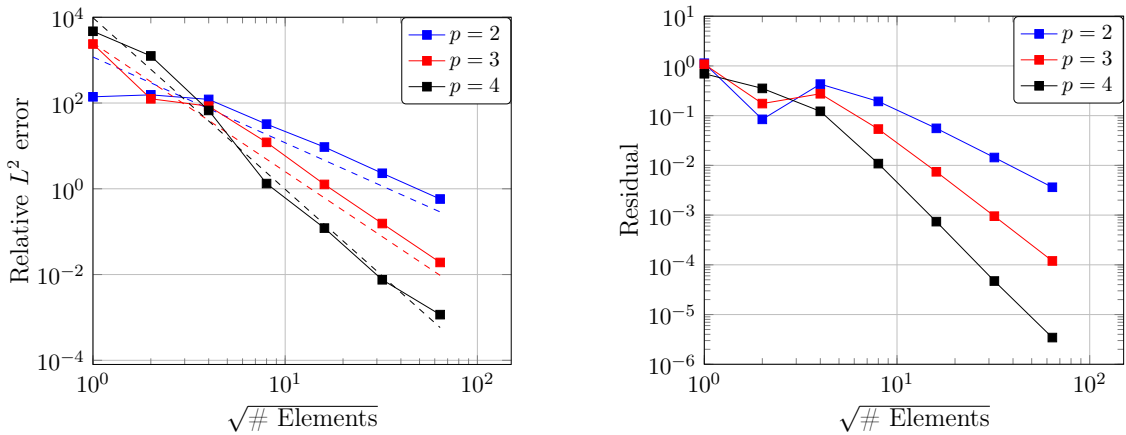
\begin{figure}[htbp]
\centering
\resizebox{0.4\textwidth}{!}{
\begin{tikzpicture}
		\begin{loglogaxis}[xmin=1,xmax=150, ymin=8e-5,ymax=10000,xlabel=\large{$\sqrt{\text{\# Elements}}$},ylabel=\large{$\text{Relative } L^2 \text{ error}$},grid=major,legend style={at={(1,1)},anchor=north east,font=\small,rounded corners=2pt}]
		\addplot[color = blue,mark=square*] table[x= sqrt_n, y=rel_L2_error, col sep = comma] {results/real_FOS/conv_p1_FOS_ellp_10p4_new_para_real.dat};
  		\addplot[dashed,color = blue,forget plot] table[x= sqrt_n, y=extslp, col sep = comma] {results/real_FOS/conv_p1_FOS_ellp_10p4_new_para_real.dat};
		\addplot [color = red,mark=square*] table[x= sqrt_n, y=rel_L2_error, col sep = comma] {results/real_FOS/conv_p2_FOS_ellp_10p4_new_para_real.dat};
    		\addplot[dashed,color = red,forget plot] table[x= sqrt_n, y=extslp, col sep = comma] {results/real_FOS/conv_p2_FOS_ellp_10p4_new_para_real.dat};
		\addplot [color = black,mark=square*] table[x= sqrt_n, y=rel_L2_error, col sep = comma] {results/real_FOS/conv_p3_FOS_ellp_10p4_new_para_real.dat};
    		\addplot[dashed,color = black,forget plot] table[x= sqrt_n, y=extslp, col sep = comma] {results/real_FOS/conv_p3_FOS_ellp_10p4_new_para_real.dat};
		\legend{$p = 2$,$p = 3$,$p = 4$}
		\end{loglogaxis}
	\end{tikzpicture}
}
\hspace{0.05\textwidth}
\resizebox{0.4\textwidth}{!}{
\begin{tikzpicture}
		\begin{loglogaxis}[xmin=1,xmax=150, ymin=1e-6,ymax=10,xlabel=\large{$\sqrt{\text{\# Elements}}$},ylabel=\large{$\text{Residual}$},grid=major,legend style={at={(1,1)},anchor=north east,font=\small,rounded corners=2pt} ]
		\addplot[color = blue,mark=square*] table[x= sqrt_n, y=res, col sep = comma] {results/real_FOS/conv_p1_FOS_ellp_10p4_new_para_real.dat};
		\addplot [color = red,mark=square*] table[x= sqrt_n, y=res, col sep = comma] {results/real_FOS/conv_p2_FOS_ellp_10p4_new_para_real.dat};
		\addplot [color = black,mark=square*] table[x= sqrt_n, y=res, col sep = comma] {results/real_FOS/conv_p3_FOS_ellp_10p4_new_para_real.dat};
		\legend{$p = 2$,$p = 3$,$p = 4$}
		\end{loglogaxis}
\end{tikzpicture}
}
\caption{First-order soliton: convergence plots for the (a) relative $L^2$ error in $\Re(u(\tau,\xi))$ and $\Re(\bm{\sigma}(\tau,\xi))$, and (b) DPG residual for the elliptic system with $c = 10000$.} \label{errConvFOSelp}
\end{figure}

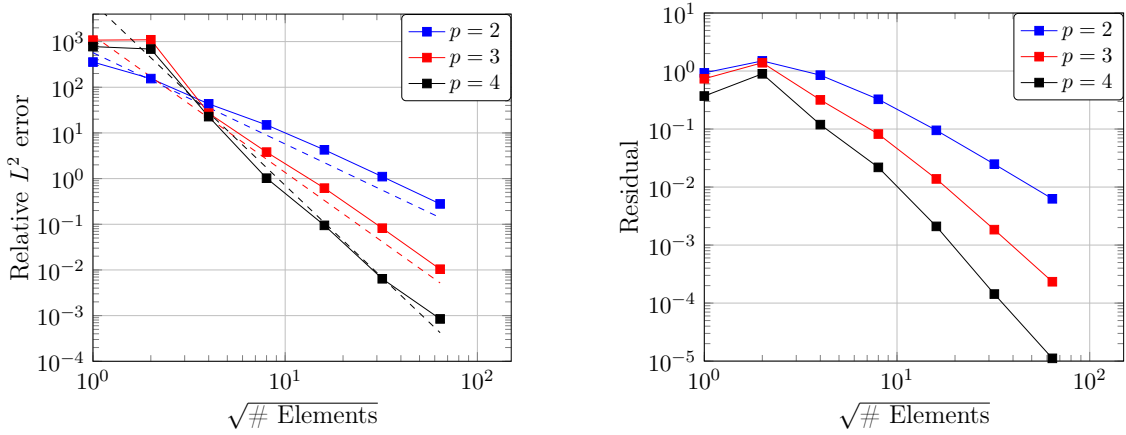
\begin{figure}[htbp]
\centering
\resizebox{0.4\textwidth}{!}{
\begin{tikzpicture}
		\begin{loglogaxis}[xmin=1,xmax=150, ymin=1e-4,ymax=4000,xlabel=\large{$\sqrt{\text{\# Elements}}$},ylabel=\large{$\text{Relative } L^2 \text{ error}$},grid=major,legend style={at={(1,1)},anchor=north east,font=\small,rounded corners=2pt}]
		\addplot[color = blue,mark=square*] table[x= sqrt_n, y=rel_L2_error, col sep = comma] {results/real_SOS/conv_p1_SOS_ellp_10p4_new_para_real.dat};
  		\addplot[dashed,color = blue,forget plot] table[x= sqrt_n, y=extslp, col sep = comma] {results/real_SOS/conv_p1_SOS_ellp_10p4_new_para_real.dat};
		\addplot [color = red,mark=square*] table[x= sqrt_n, y=rel_L2_error, col sep = comma] {results/real_SOS/conv_p2_SOS_ellp_10p4_new_para_real.dat};
    		\addplot[dashed,color = red,forget plot] table[x= sqrt_n, y=extslp, col sep = comma] {results/real_SOS/conv_p2_SOS_ellp_10p4_new_para_real.dat};
		\addplot [color = black,mark=square*] table[x= sqrt_n, y=rel_L2_error, col sep = comma] {results/real_SOS/conv_p3_SOS_ellp_10p4_new_para_real.dat};
    		\addplot[dashed,color = black,forget plot] table[x= sqrt_n, y=extslp, col sep = comma] {results/real_SOS/conv_p3_SOS_ellp_10p4_new_para_real.dat};
		\legend{$p = 2$,$p = 3$,$p = 4$}
		\end{loglogaxis}
	\end{tikzpicture}
}
\hspace{0.05\textwidth}
\resizebox{0.4\textwidth}{!}{
\begin{tikzpicture}
		\begin{loglogaxis}[xmin=1,xmax=150, ymin=1e-5,ymax=10,xlabel=\large{$\sqrt{\text{\# Elements}}$},ylabel=\large{$\text{Residual}$},grid=major,legend style={at={(1,1)},anchor=north east,font=\small,rounded corners=2pt} ]
		\addplot[color = blue,mark=square*] table[x= sqrt_n, y=res, col sep = comma] {results/real_SOS/conv_p1_SOS_ellp_10p4_new_para_real.dat};
		\addplot [color = red,mark=square*] table[x= sqrt_n, y=res, col sep = comma] {results/real_SOS/conv_p2_SOS_ellp_10p4_new_para_real.dat};
		\addplot [color = black,mark=square*] table[x= sqrt_n, y=res, col sep = comma] {results/real_SOS/conv_p3_SOS_ellp_10p4_new_para_real.dat};
		\legend{$p = 2$,$p = 3$,$p = 4$}
		\end{loglogaxis}
\end{tikzpicture}
}
\caption{Second-order soliton: convergence plots for the (a) relative $L^2$ error in $\Re(u(\tau,\xi))$ and $\Re(\bm{\sigma}(\tau,\xi))$, and (b) DPG residual for the elliptic system with $c = 10000$.} \label{errConvSOSelp}
\end{figure}

Subsequently, these systems are solved computationally under various element discretization sizes ($h$), and for polynomial orders $p = 2, 3, 4$, and are compared against their respective analytic solutions in order to extract the numerical convergence rates. 
The observed convergence rates for the hyperbolic system are shown in~\Cref{errConvFOShyp,errConvSOShyp}, and, for the elliptic system, in~\Cref{errConvFOSelp,errConvSOSelp}.
These results effectively demonstrate that our approach does achieve the expected optimal convergence rates $(p+1)$ for both soliton orders, in both the hyperbolic and elliptic systems. 

\begin{figure}[htbp]
\centering
\resizebox{0.4\textwidth}{!}{
\begin{tikzpicture}
		\begin{loglogaxis}[xmin=1,xmax=150, ymin=1e-3,ymax=200,xlabel=\large{$\sqrt{\text{\# Elements}}$},ylabel=\large{$\text{Relative } L^2 \text{ error}$},grid=major,legend style={at={(1,1)},anchor=north east,font=\small,rounded corners=2pt}]
		\addplot[color = blue,mark=square*] table[x= sqrt_n, y=rel_L2_error, col sep = comma] {results/conv_p1_FOS_hyp_new_para.dat};
  		\addplot[dashed,color = blue,forget plot] table[x= sqrt_n, y=extslp, col sep = comma] {results/conv_p1_FOS_hyp_new_para.dat};
		\addplot [color = red,mark=square*] table[x= sqrt_n, y=rel_L2_error, col sep = comma] {results/conv_p2_FOS_hyp_new_para.dat};
    		\addplot[dashed,color = red,forget plot] table[x= sqrt_n, y=extslp, col sep = comma] {results/conv_p2_FOS_hyp_new_para.dat};
		\addplot [color = black,mark=square*] table[x= sqrt_n, y=rel_L2_error, col sep = comma] {results/conv_p3_FOS_hyp_new_para.dat};
    		\addplot[dashed,color = black,forget plot] table[x= sqrt_n, y=extslp, col sep = comma] {results/conv_p3_FOS_hyp_new_para.dat};
		\legend{$p = 2$,$p = 3$,$p = 4$}
		\end{loglogaxis}
	\end{tikzpicture}
}
\hspace{0.05\textwidth}
\resizebox{0.4\textwidth}{!}{
\begin{tikzpicture}
		\begin{loglogaxis}[xmin=1,xmax=150, ymin=5e-5,ymax=1,xlabel=\large{$\sqrt{\text{ \# Elements}}$},ylabel=\large{$\text{Residual}$},grid=major,legend style={at={(1,1)},anchor=north east,font=\small,rounded corners=2pt} ]
		\addplot[color = blue,mark=square*] table[x= sqrt_n, y=res, col sep = comma] {results/conv_p1_FOS_hyp_new_para.dat};
		\addplot [color = red,mark=square*] table[x= sqrt_n, y=res, col sep = comma] {results/conv_p2_FOS_hyp_new_para.dat};
		\addplot [color = black,mark=square*] table[x= sqrt_n, y=res, col sep = comma] {results/conv_p3_FOS_hyp_new_para.dat};
		\legend{$p = 2$,$p = 3$,$p = 4$}
		\end{loglogaxis}
\end{tikzpicture}
}
\caption{First-order soliton: convergence plots for the (a) relative $L^2$ error in $\Re(u(\tau,\xi))$ and $\Re(v(\tau,\xi))$, and (b) DPG residual for the hyperbolic system.} \label{errConvFOShyp}
\end{figure}
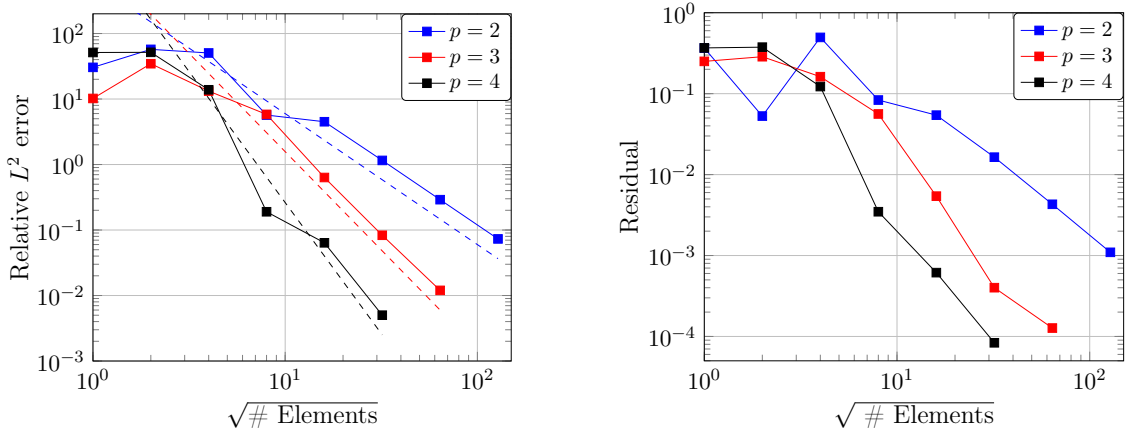

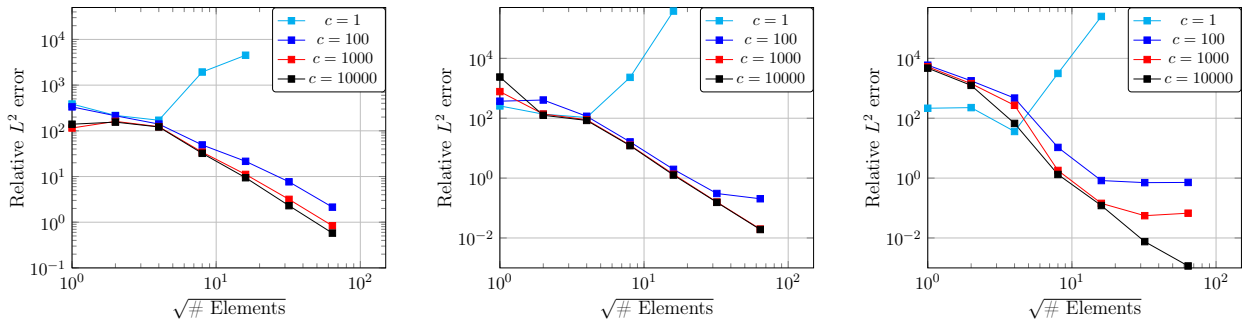
\begin{figure}[htbp]
\centering
\resizebox{0.3\textwidth}{!}{
\begin{tikzpicture}
		\begin{loglogaxis}[xmin=1,xmax=150, ymin=1e-1,ymax=50000,xlabel=\large{$\sqrt{\text{\# Elements}}$},ylabel=\large{$\text{Relative } L^2 \text{ error}$},grid=major,legend style={at={(1,1)},anchor=north east,font=\small,rounded corners=2pt}]
  		\addplot[color = cyan,mark=square*] table[x= sqrt_n, y=rel_L2_error, col sep = comma] {results/real_FOS/conv_p1_FOS_ellp_10p0_new_para_real.dat};
		\addplot[color = blue,mark=square*] table[x= sqrt_n, y=rel_L2_error, col sep = comma] {results/real_FOS/conv_p1_FOS_ellp_10p2_new_para_real.dat};
		\addplot [color = red,mark=square*] table[x= sqrt_n, y=rel_L2_error, col sep = comma] {results/real_FOS/conv_p1_FOS_ellp_10p3_new_para_real.dat};
		\addplot [color = black,mark=square*] table[x= sqrt_n, y=rel_L2_error, col sep = comma] {results/real_FOS/conv_p1_FOS_ellp_10p4_new_para_real.dat};
		\legend{$c=1$,$c = 100$,$c = 1000$,$c = 10000$}
		\end{loglogaxis}
	\end{tikzpicture}
}
\hspace{0.01\textwidth}
\resizebox{0.3\textwidth}{!}{
\begin{tikzpicture}
		\begin{loglogaxis}[xmin=1,xmax=150, ymin=1e-3,ymax=500000,xlabel=\large{$\sqrt{\text{\# Elements}}$},ylabel=\large{$\text{Relative } L^2 \text{ error}$},grid=major,legend style={at={(1,1)},anchor=north east,font=\small,rounded corners=2pt}]
  		\addplot[color = cyan,mark=square*] table[x= sqrt_n, y=rel_L2_error, col sep = comma] {results/real_FOS/conv_p2_FOS_ellp_10p0_new_para_real.dat};
		\addplot[color = blue,mark=square*] table[x= sqrt_n, y=rel_L2_error, col sep = comma] {results/real_FOS/conv_p2_FOS_ellp_10p2_new_para_real.dat};
		\addplot [color = red,mark=square*] table[x= sqrt_n, y=rel_L2_error, col sep = comma] {results/real_FOS/conv_p2_FOS_ellp_10p3_new_para_real.dat};
		\addplot [color = black,mark=square*] table[x= sqrt_n, y=rel_L2_error, col sep = comma] {results/real_FOS/conv_p2_FOS_ellp_10p4_new_para_real.dat};
		\legend{$c=1$,$c = 100$,$c = 1000$,$c = 10000$}
		\end{loglogaxis}
	\end{tikzpicture}
}
\hspace{0.01\textwidth}
\resizebox{0.3\textwidth}{!}{
\begin{tikzpicture}
		\begin{loglogaxis}[xmin=1,xmax=150, ymin=1e-3,ymax=500000,xlabel=\large{$\sqrt{\text{\# Elements}}$},ylabel=\large{$\text{Relative } L^2 \text{ error}$},grid=major,legend style={at={(1,1)},anchor=north east,font=\small,rounded corners=2pt}]
  		\addplot[color = cyan,mark=square*] table[x= sqrt_n, y=rel_L2_error, col sep = comma] {results/real_FOS/conv_p3_FOS_ellp_10p0_new_para_real.dat};
		\addplot[color = blue,mark=square*] table[x= sqrt_n, y=rel_L2_error, col sep = comma] {results/real_FOS/conv_p3_FOS_ellp_10p2_new_para_real.dat};
		\addplot [color = red,mark=square*] table[x= sqrt_n, y=rel_L2_error, col sep = comma] {results/real_FOS/conv_p3_FOS_ellp_10p3_new_para_real.dat};
		\addplot [color = black,mark=square*] table[x= sqrt_n, y=rel_L2_error, col sep = comma] {results/real_FOS/conv_p3_FOS_ellp_10p4_new_para_real.dat};
		\legend{$c=1$,$c = 100$,$c = 1000$,$c = 10000$}
		\end{loglogaxis}
	\end{tikzpicture}
}
\caption{Second-order soliton: comparison of convergence for the elliptic system of equations for different values of scaling constant $c$ using (a) $p =1$, (b) $p =2$, and $p =3$.} \label{err_scaling_comp}
\end{figure}

As another demonstration of expected model behavior, the elliptic system is re-computed for different values of the scaling constant $c$, at various polynomial orders, to emphasize that proper conditioning is needed for convergence. 
As noted in~\cref{sec:DPG_disc}, for $c = 1$, we obtain a highly ill-conditioned system because of the poor scaling of the elliptic system equations. 
However, as $c$ is increased, improved conditioning in the linear system is observed, leading to the convergence behavior illustrated in~\Cref{err_scaling_comp}. 

\begin{figure}[htbp]
\centering
\resizebox{0.4\textwidth}{!}{
\begin{tikzpicture}
		\begin{loglogaxis}[xmin = 10, xmax = 1000, ymin = 1e-2, ymax = 5000, xlabel = \large{$\sqrt{\text{ndof}}$}, ylabel = \large{$\text{Relative } L^2 \text{ error}$}, grid = major, legend style = {at={(1,1)},anchor=north east,font=\small,rounded corners=2pt}]
		\addplot[color = blue,mark=square*] table[x= sqrt_ndof, y=rel_L2_error, col sep = comma] {results/error_uniform_beam.dat};
  
		\addplot [color = red,mark=square*] table[x= sqrt_ndof, y=rel_L2_error, col sep = comma] {results/error_adaptation_beam.dat};

		\legend{$p = 3 (\text{uniform })$,$p = 3 (\text{adaptive})$}
		\end{loglogaxis}
	\end{tikzpicture}
}
\hspace{0.05\textwidth}
\resizebox{0.4\textwidth}{!}{
\begin{tikzpicture}
		\begin{loglogaxis}[xmin=10,xmax=1000, ymin=5e-5,ymax=1,xlabel=\large{$\sqrt{\text{ndof}}$},ylabel=\large{$\text{Residual}$},grid=major,legend style={at={(1,1)},anchor=north east,font=\small,rounded corners=2pt} ]
		\addplot[color = blue,mark=square*] table[x= sqrt_ndof, y=res, col sep = comma] {results/error_uniform_beam.dat};
  
		\addplot [color = red,mark=square*] table[x= sqrt_ndof, y=res, col sep = comma] {results/error_adaptation_beam.dat};

		\legend{$p = 3 (\text{uniform})$,$p = 3 (\text{adaptive})$}
		\end{loglogaxis}
\end{tikzpicture}
}
\caption{Adaptive refinements: convergence plots for the (a) relative $L^2$ error in $\Re(u(\tau,\xi))$ and $\Re(\bm{\sigma}(\tau,\xi))$, and (b) DPG residual for the elliptic system with $c = 10000$.} \label{errConvSoladap}
\end{figure}
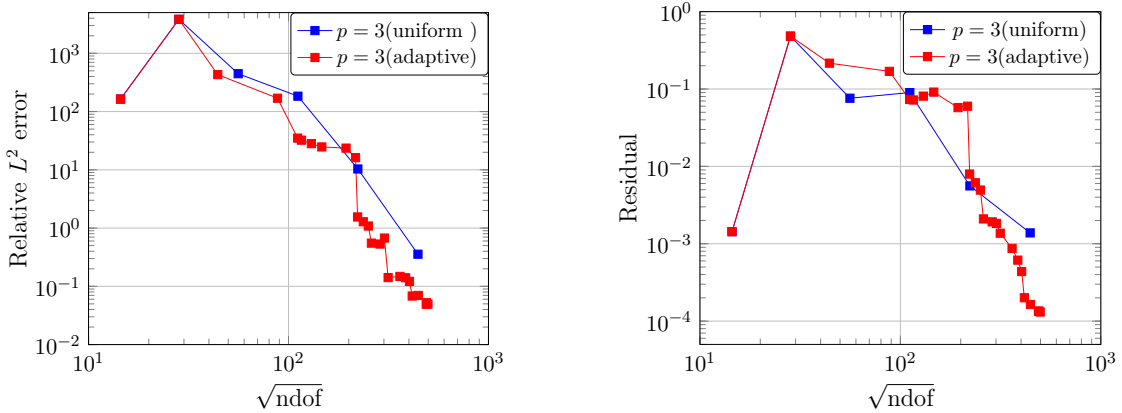

Finally, in~\Cref{errConvSoladap}, we demonstrate the capability of the DPG methodology to drive automatic mesh adaptivity using its built-in residual-based error estimator. Here, the corresponding analytical solution for the model that substantiates this adaptivity is given by 
\begin{equation}
     u(\tau,\xi) = \frac{1}{\sqrt{2\pi \omega}}e^{-\frac{1}{2}\frac{{(\tau-0.5)}^2}{\omega}} + i \sin(\xi),
\end{equation}
where $\omega = 0.001$. 

\section{Conclusions} \label{sec: conclusions}

In this article, we have proposed and validated a new space-time finite element model for the pulsed laser propagation in optical waveguides. 
For the sake of well-posedness, our model extends the canonical NLS equation by retaining the second-order spatial derivative of the amplitude; presenting it as a small perturbation rather than neglecting it under the slowly-varying envelope approximation (or paraxial-wave approximation), as is often done in the optics community. 
This model is then re-posed as a stable first-order system of hyperbolic or elliptic equations, depending upon the sign of the group velocity dispersion ($\beta_2$). 
We prove the stability and well-posedness of both the hyperbolic and elliptic systems, and further support this through various numerical examples. 
While this article focuses on the well-posedness of our model in terms of the resulting differential operator, further analysis is needed to address the nonlinear load. 
This aspect will be explored in future work. 
\subsection*{Declarations}
The authors have no relevant financial or non-financial interests to disclose.
\paragraph{Disclaimers.}
\addcontentsline{toc}{section}{Disclaimers}
This article has been approved for public release; distribution unlimited. Public Affairs release approval {\#}AFRL-2024-6442. 
The views expressed in this article are those of the authors and do not necessarily reflect the official policy or position of the Department of the Air Force, the Department of Defense, or the U.S. government.


\printbibliography[heading=bibintoc]
\setcounter{section}{0}
\end{document}